\documentclass[smallextended,referee,envcountsect,]{svjour3}
\smartqed
\usepackage{graphicx}
\usepackage{subfig}
\journalname{JOTA}

\usepackage{calrsfs}
\DeclareMathAlphabet{\pazocal}{OMS}{zplm}{m}{n}

\usepackage{enumitem}
\usepackage[utf8]{inputenc} 
\usepackage[T1]{fontenc}    
\usepackage{graphicx}       
\usepackage{hyperref}       
\usepackage{amsmath}        
\usepackage{amssymb}        
\usepackage{fancyhdr}       
\usepackage{lastpage}       
\usepackage{geometry}       
\usepackage{cite} 
\usepackage[dvipsnames]{xcolor}
\usepackage[english]{babel}
\usepackage{comment}
\newtheorem{assumption}{Assumption}
\newcommand{\col}{\text{col}}
\usepackage{cancel}
\usepackage{algorithm, algorithmic}
\newtheorem{thrm}{Theorem}
\DeclareMathOperator*{\argmin}{arg\,min}
\usepackage{xcolor}
\newcommand{\GP}[1]{\textcolor{black}{#1}}
\newcommand{\SG}[1]{\textcolor{black}{#1}}
\newcommand{\RB}[1]{\textcolor{black}{#1}}

\newcommand{\x}{\bar x}
\newcommand{\n}[1]{\|#1 \|}

\begin{document}

\title{ Nash equilibrium seeking for a class of quadratic-bilinear Wasserstein distributionally robust games}
\author{\makebox{Georgios Pantazis, Reza Rahimi Baghbadorani, Sergio Grammatico}\thanks{The authors would like to thank Prof. Dimitris Boskos for his useful feedback and discussions during the preparation of this work.}}
\institute{Delft Center for Systems and Control (DCSC)
  \{G.Pantazis, R.Rahimibaghbadorani,S.Grammatico\}@tudelft.nl.
}

\date{Received: date / Accepted: date}

\maketitle
\begin{abstract}
	We consider a class of Wasserstein distributionally robust Nash equilibrium problems, where agents construct heterogeneous data-driven Wasserstein ambiguity sets using private samples and radii, in line with their individual risk-averse preference. By leveraging relevant properties of this class of games, we show that equilibria of the original seemingly infinite-dimensional problem can be obtained as a solution to a finite-dimensional Nash equilibrium problem. We then reformulate the associated variational problem as a finite-dimensional variational inequality and establish the connection between the corresponding solution sets. 
 Our reformulation has scalable behaviour with respect to the data size and maintains a fixed number of constraints, independently of the number of samples.
 To numerically compute a solution, we leverage two algorithms, based on the golden ratio algorithm. 
 The efficiency of both algorithmic schemes is corroborated through extensive simulation studies on an illustrative example and a stochastic portfolio allocation game, where behavioural coupling among investors is modeled. 
\end{abstract}


\section{Introduction}
A wide range of applications, from smart grids \cite{Saad1} and communication networks \cite{Scutari} to social networks \cite{Acemoglu2013} can be modeled as a collection of self-interested interacting decision makers optimizing different objective functions under operational constraints. Game theory \cite{Basar1} provides the fundamental theoretical framework for analyzing such systems. Although investigating deterministic games can be adequate in some case studies \cite{Scutari}, \cite{Paccagnan2017}, most real-world applications involve decision making under uncertainty, which stresses the need for the inclusion of stochasticity in the existing models. Several studies have explored uncertainty within a game theoretic context, based on particular technical assumptions on the probability distribution \cite{Kouvaritakis}, \cite{Singh} and/or the properties of the uncertainty sample space \cite{Aghassi2006, FukuSOCCP}.  \par 

When the probability distribution is unknown and distribution models are not an accurate description of the stochastic aspect of the problem, sampling-based or data-driven methods have shown strong potential for proposing robust solutions against uncertainty. Works such as \cite{Feleconf2019, Fele2021, Dario_Scenario, fele-a, Pantazis2020, mammarela2023, Pantazis2023_apriori} design distribution-free approaches for data-driven Nash equilibria based on statistical learning techniques. More specifically, \cite{fele-a, Pantazis2020, mammarela2023, Pantazis2023_apriori} account for possible strategic perturbations around the Nash equilibrium. Separately, works based on Sample Average Approximation (SAA) techniques, such as \cite{Franci_2021, Franci_2021_merely}, develop algorithms for finding Nash equilibria in stochastic settings by using expected values of cost functions. The works mentioned above constitute data-driven methods for stochastic equilibrium seeking. These works, however, do not account for ambiguity in the probability distribution, where the distribution itself may be uncertain within some known bounds. The challenge of ambiguity in the distributions becomes pronounced in multi-agent settings, where heterogeneous uncertainties affect the agents' costs, often necessitating the consideration of different ambiguity sets, each representing the individual risk-averse nature of each agent. 
\par To account for distributional uncertainty, distributionally robust optimization (DRO) uses a so-called ambiguity set of possible probability distributions to make decisions robust against probabilistic variations within this set. Unlike scenario-based methods, which require many samples for robustness, DRO can perform well with limited data by adjusting the ambiguity set. DRO includes special cases like sample average approximation (SAA) and robust optimization (RO). At the same time, DRO can be less conservative than RO and \SG{offer} better out-of-sample performance than SAA, making it especially useful in data-driven applications with limited data. Recently, Wasserstein ambiguity sets \cite{villani_topics_2016}, which use empirical data and the Wasserstein metric to measure distributional deviations, have gained attention. These sets are favored for penalizing horizontal shifts and providing finite-sample guarantees. Research has focused on convergence of empirical estimates in the Wasserstein distance \cite{Dereich, mohajerin_esfahani_data-driven_2018, Dedecker1, Weed,Weed_2, Fournier_2023}, \SG{as well as obtaining} tractable reformulations of Wasserstein distributionally robust optimization problems \cite{mohajerin_esfahani_data-driven_2018, netessine_wasserstein_2019, Lotfi1, Lotfi2}.   Extensions \SG{of those works} include distributionally robust chance-constrained programs \cite{Chen2018, Hota2018, Alamo2024risk}. \par  
 Despite the considerable body of literature on DRO with Wasserstein ambiguity sets, the exploration of data-driven Wasserstein distributionally robust Nash equilibrium problems with heterogeneous uncertainty in the cost functions represents a notably underexplored topic. Most works in the literature consider moment-based methods or other measures of distance between distributions.  For instance, \cite{Peng2021} considers a non-cooperative game with distributionally robust chance-constrained strategy sets applied to duopoly Cournot competition. \GP{ The work \cite{Liu2018} develops distributionally robust equilibrium models based on the Kullback-Leibler (KL) divergence for hierarchical competition in supply chains.} Other works mainly consider ambiguity in the constraints,  such as the recent work \cite{Xia_elliptical_2023}, which studies a game with deterministic cost for each agent and  distributionally robust chance constraints with the \RB{centre} of the Wasserstein ambiguity set being an elliptical distribution;  \cite{fabiani2023distributionally} reformulates an  equilibrium problem with a deterministic cost and distributionally robust chance-constraints as a mixed-integer generalized Nash equilibrium problem leveraging the results in \cite{Chen2023}. 
The contributions of this work with respect to the related literature are the following:
\begin{enumerate}[label=(\roman*)]
\item \RB{We study a class of heterogeneous data-driven Wasserstein distributionally robust games with a quadratic-bilinear structure for each agent objective function. Each agent is considered to act as a selfish entity making individual decisions to optimize  their own cost function, which depends on a private uncertain parameter. The game-theoretic formulation generalizes the optimization framework to allow for interdependence of the agent cost functions. Furthermore, as we model possibly heterogeneous probability distributions, we lift the assumption of a common ambiguity set, allowing  each agent to construct their own, based on their private data and personalized Wasserstein radius. In this generalized setting, the solution set can significantly differ from those associated with common ambiguity sets \cite{mohajerin_esfahani_data-driven_2018}, \cite{netessine_wasserstein_2019}, \cite{shafiee2024}, \cite{Lanzetti2025}  } 

\item We reformulate the original game as a robust Nash equilibrium problem and establish the connection between the distributionally robust and robust Nash equilibria of the corresponding problems.  For this class of games, we demonstrate that the inner maximization can be solved without the use of epigraphic variables \cite{kuhn2019wasserstein}, \cite{pantazis2023_DRG}, which, unlike distributionally robust optimization, in distributionally robust games can lead to unshared coupling constraints. As such, our approach decreases computational complexity significantly. To the best of our knowledge, this is the first distributionally robust game-theoretic reformulation that leads to data-scalable results by leveraging the structure of the problem at hand. 
\item The robust Nash equilibrium problem is then reformulated as a variational inequality (VI). Unlike results of similar classes of problems in optimization \cite{Boskos_2024}, where the reformulated variational inequality is monotone under certain assumptions, the mapping corresponding to the game can be nonmonotone in general due to the heterogeneity of the agents' ambiguity sets and costs.  However, we show that this problem can be efficiently solved empirically using two algorithms: the adaptive golden ratio algorithm (aGRAAL) \cite{malitsky_golden_2020} and a hybrid version of this algorithm ({Hybrid-Alg} in \cite{Reza_2024}). Notably, our numerical results show that in several cases, the convergence speed is close to linear, and increasing the number of samples does not slow down the convergence. Our results are then applied to a portfolio allocation game that takes into account market uncertainty and behavioural coupling of market participants. 
\end{enumerate}

 \section{Problem Formulation}
 \emph{Notation}:
In this section, we introduce some basic notation and results required for the subsequent developments. To this end, consider the index set  $\mathcal{N}=\{1, \dots, N\}$. The decision vector of each agent $i \in \mathcal{N}$ is denoted by $x_i=\col((x^{(j)}_i)_{j=1}^{n}) \in X_i \subseteq \mathbb{R}^{n}$, where $x^{(j)}_i, j=1,\dots,n$, denotes an element of the decision vector; let $x_{-i}=\col((x_j)_{j=1, j \neq i}^N) \in X_{-i}=\prod_{i=1, i\neq j}^N X_j \subseteq \mathbb{R}^{(N-1)n}$ be the decision vector of all other agents' decisions except for that of agent $i$ and $x=\col((x_i)_{i =1}^N)$ be the collective decision vector.  We denote $||\cdot||=\|\cdot\|_2$.  The projection operator $\text{proj}_X(x)$ of a point $x$ to the set $X$ is given by $\text{proj}_X(x)= \argmin_{y \in X}\|x-y\|$.   	$F$ is monotone on $X$ if $(x-y)^\top(F(x)-F(y)) \geq 0 $ for all $x, y \in X$.   If the condition is not satisfied, the mapping is called nonmonotone.
 Let us  denote $\pazocal{P}(\mathbb{R}^m)$ as the set of all probability measures on $\mathbb{R}^m$ and  define
\begin{align*}
\pazocal{M}(\mathbb{R}^m)=\left\{ \mathbb{Q} \mid \mathbb{Q} \text{ is a distribution on } \mathbb{R}^m  \text{ and }  \mathbb{E}_\mathbb{P}[ \|\xi\| ]=\int_\Xi \|\xi \| \mathbb{Q}(d\xi) < \infty \right\} . 
 \end{align*}
In other words, $\pazocal{M}(\mathbb{R}^m)$ considers the sets of all distributions defined on $\mathbb{R}^m$ with a bounded first-order moment. 
We are now ready to  define the Wasserstein metric to quantify the distance between two probability distributions.


\begin{definition} \label{Wasserstein}
	The Wasserstein metric $d_W: \pazocal{M}(\mathbb{R}^m) \times \pazocal{M}(\mathbb{R}^m) \rightarrow \mathbb{R}_{\geq 0}$  between two distributions $\mathbb{Q}_1, \mathbb{Q}_2 \in \pazocal{M}(\mathbb{R}^m)$ is defined as 
	\begin{align}
		d_W(\mathbb{Q}_1, \mathbb{Q}_2):=&\left (\inf_{\Pi \in \pazocal{J}(\xi_1 \sim \mathbb{Q}_1, \xi_2 \sim \mathbb{Q}_2)}\int_{\mathbb{R}^m \times \mathbb{R}^m } \|\xi_1-\xi_2 \|^p \Pi(d\xi_1, d\xi_2) \right )^{1/p} , \nonumber 
	\end{align}
where  $\pazocal{J}(\xi_1 \sim \mathbb{Q}_1, \xi_2 \sim \mathbb{Q}_2)$ represent the set of joint probability distributions of the random variables $\xi_1$ and $\xi_2$ with marginals  $\mathbb{Q}_1$ and $\mathbb{Q}_2$, respectively.  \hfill $\square$
\end{definition}
The Wasserstein metric can be viewed as the optimal transport plan to fit the probability distribution $\mathbb{Q}_1$ to $\mathbb{Q}_2$ \cite{villani_topics_2016}. 
\subsection{Problem Formulation}
Consider a population of agents with index set $\mathcal{N}=\{1, \dots, N\}$. Based on that we define the following game:

\begin{align}
\forall i \in \mathcal{N}:	 \min_{x_i \in X_i} \max_{\mathbb{Q}_i \in \mathcal{P}_i} \{f_i(x_i, x_{-i})+ \mathbb{E}_{\xi_i \sim \mathbb{Q}_i}[g_i(x_i, x_{-i}, \xi_i)] \}, \nonumber 
\end{align}
where $f_i: \mathbb{R}^{n} \times \mathbb{R}^{n(N-1)}  \rightarrow \mathbb{R} $, $g_i:  \mathbb{R}^{n} \times \mathbb{R}^{n(N-1)} \times \mathbb{R}^m  \rightarrow \mathbb{R}$  for all $i \in \mathcal{N}$ and $\mathcal{P}_i$ is the ambiguity set of the uncertain parameter $\xi_i$.  We call the collection of the coupled optimization problems above for all agents $i \in \mathcal{N}$ as game ($G$).

For game $G$, we define the notion of \emph{distributionally robust Nash equilibrium} as follows:
\begin{definition}
	A decision vector $x^\ast \in \prod_{i=1}^N X_i$ is a distributionally robust Nash equilibrium (DRNE) of game $G$ if, given the decisions of all other agents $x^\ast_{-i} \in X_{-i}$ it holds that
	\begin{align} \label{DR_formulation}
		x^\ast_i \in \argmin_{x_i \in X_i } \max_{\mathbb{Q}_i \in \mathcal{P}_i} \{  f_i(x_i, x^\ast_{-i})+ \mathbb{E}_{\xi_i \sim \mathbb{Q}_i}[g_i(x_i, x^\ast_{-i}, \xi_i)] \}, \ \forall i \in \mathcal{N}.
	\end{align}
\end{definition}

In other words, a decision vector \( x^\ast \) is a DRNE of $G$ if, for each agent \( i \), given the equilibrium strategies \( x^\ast_{-i} \) of all other agents with their respective local sets, the following holds: player \( i \) chooses their strategy \( x_i \) in a way that minimizes their objective, considering both their deterministic cost $f_i$ and the worst-case expected effect of distributional uncertainty in $\xi_i$ of \( g_i \). This must hold for all agents simultaneously, ensuring that no agent can improve their outcome by unilaterally changing their strategy, even in the face of worst-case distribution. 

In this work, we follow a data-driven approach and consider \emph{heterogeneous} Wasserstein ambiguity sets constructed by each individual agent on  the basis of their own individual data. Thus, we define an appropriate notion of distance between probability distributions. Due to its ability to penalize horizontal dislocations of distributions and often capturing realistic shifts in distributions, in this work we will use the Wasserstein distance. Specifically,  for each $i \in \mathcal{N}$ the empirical probability distribution is constructed based on $K_i$ independent and identically distributed (i.i.d.) samples $\xi_{K_i}=\{\xi^{(1)}_i,  \dots, \xi^{(K_1)}_i\}$ drawn by agent $i$ as follows: 
	\[
	\hat{\mathbb{P}}_{K_i} = \frac{1}{K_i} \sum_{k_i=1}^{K_i} \delta_{\xi^{(k_i)}_i}
	\]
	where \(\delta_{\xi_i}\) is the Dirac delta measure that assigns the full probability mass at the point \(\xi_i\). We then consider a radius $\varepsilon_i$, based on the Wasserstein distance, and construct the data-driven Wasserstein ambiguity ball of agent $i$ as follows:
	\begin{align} \label{ambiguity_set}
		\mathcal{P}_i=\{ Q_i \in \mathcal{P}(\mathbb{R}^m) \mid d_W(Q_i , \hat{\mathbb{P}}_{K_i} ) \leq \varepsilon_i\},
	\end{align}
where $\mathcal{P}(\mathbb{R}^m)$ denotes the collection of all probability distributions defined on the support set $\mathbb{R}^m$. 

\subsection{\RB{Effect of heterogeneous  ambiguity sets on the equilibrium set}}

We consider a simple example to illustrate how heterogeneity of each agents' ambiguity sets, which correspond to an uncertain parameter $\xi \in \mathbb{R}^2$, can affect the solution set of DRNE even in the case of 1-Wasserstein distance $(p=1)$. In the following distributionally robust game,  two self-interested agents solve the following interdependent optimization problems:

\begin{align}
\begin{cases}
 \min\limits_{x_1 }\max\limits_{\mathbb{Q}_1 \in \mathcal{P}_{\varepsilon_1}(\hat{\mathbb{P}}_{K_1})} \mathbb{E}_{\mathbb{Q}_1}[c_{11}x^2_1+c_{12}x_1x_2+\xi^\top x]  \\
\min\limits_{x_2 }\max\limits_{\mathbb{Q}_2 \in \mathcal{P}_{\varepsilon_2}(\hat{\mathbb{P}}_{K_2})} \mathbb{E}_{\mathbb{Q}_2}[c_{21}x_1x_2+c_{22}x^2_2+\xi^\top x] \end{cases} \label{example}
\end{align}
The game in (\ref{example}) can be equivalently written as follows \cite{netessine_wasserstein_2019}:
\begin{align}
\begin{cases}
 \min\limits_{x_1} c_{11}x_1^2 + c_{12}x_1x_2 + \frac{1}{K_1} \sum\limits_{k=1}^{K_1} (\xi^{(k)}_1)^\top x + \varepsilon_1  \|x\|_* \nonumber \\
 \min\limits_{x_2} c_{21}x_1x_2 + c_{22}x_2^2 + \frac{1}{K_2} \sum\limits_{k=1}^{K_2} (\xi^{(k)}_2)^\top x + \varepsilon_2  \|x\|_*, 
\end{cases}
\end{align}
where $\|\cdot\|_*$ represents the dual norm of $\|\cdot\|$.
At an equilibrium $\bar{x}=(\bar{x}_1, \bar{x}_2)$, the following stationarity conditions hold:
\begin{align}
\begin{cases}
2c_{11}\bar{x}_1+c_{12}\bar{x}_2+\varepsilon_1  \partial_{x_1} \|\bar{x}\|_{*}=-\frac{1}{K_1}\sum\limits_{k=1}^{K_1}\xi^{(k)}_1 \\
c_{21}\bar{x}_1+2c_{22}\bar{x}_2+\varepsilon_2  \partial_{x_2} \|\bar{x}\|_{*}=-\frac{1}{K_2}\sum\limits_{k=1}^{K_2}\xi^{(k)}_2
 \end{cases} \label{conditions}
\end{align}

Considering the 2-norm for the transport cost, we have that:

\begin{figure}[t]
\centering
    \includegraphics[scale=0.45]{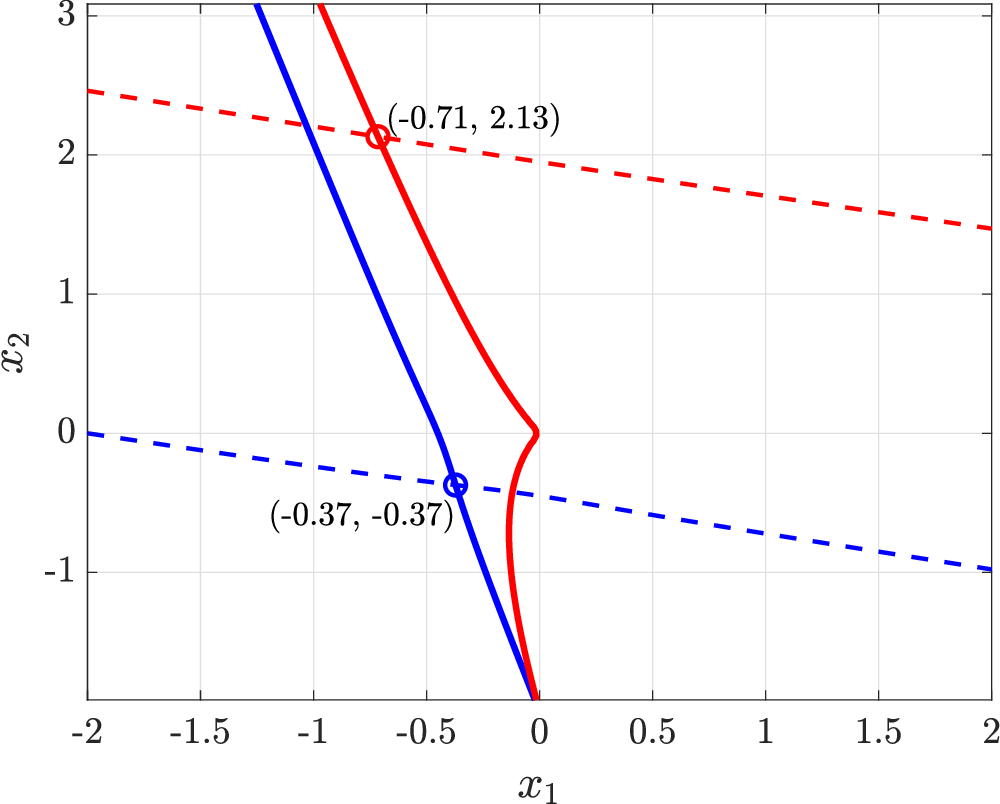}
    \hspace{2mm}
    \caption{Equilibrium solution for homogeneous (blue lines) and heterogeneous (red lines) ambiguity sets for a 2-player game with scalar decisions. The dashed lines correspond to the second equation in (5). } \label{solution_set}
\end{figure}

\begin{align}
\begin{cases}
2c_{11}\bar{x}_1+c_{12}\bar{x}_2+\varepsilon_1 \frac{\bar{x}_1} {\sqrt{\bar{x}_1^2+\bar{x}_2^2}}=-\frac{1}{K_1}\sum\limits_{k=1}^{K_1}\xi^{(k)}_1 \\
c_{21}\bar{x}_1+2c_{22}\bar{x}_2+\varepsilon_2 \frac{\bar{x}_2}{\sqrt{\bar{x}_1^2+\bar{x}_2^2}}=-\frac{1}{K_2}\sum\limits_{k=1}^{K_2}\xi^{(k)}_2
 \end{cases} \label{conditions}
\end{align}
We now show how, even for simple discrete distributions concentrated at one point, a difference between $\mathbb{P}_1$ and $\mathbb{P}_2$ can lead to a significant change to the equilibrium solution. Consider a  discrete probability distribution $\mathbb{P}_1$ with full probability mass at $x=p_1$ (i.e., with probability 1) and 0 elsewhere, and the probability distribution $\mathbb{P}_2$ with full probability mass at $x=p_2$, 0 elsewhere:
\begin{align}
\begin{cases}
2c_{11}\bar{x}_1+c_{12}\bar{x}_2+\varepsilon_1 \frac{\bar{x}_1} {\sqrt{\bar{x}_1^2+\bar{x}_2^2}}=-p_1 \\
c_{21}\bar{x}_1+2c_{22}\bar{x}_2+\varepsilon_2 \frac{\bar{x}_2}{\sqrt{\bar{x}_1^2+\bar{x}_2^2}}=-p_2
 \end{cases} \label{conditions}
\end{align}
In Figure \ref{solution_set}, we first consider $c_{11} = 1$, $c_{12} = 0.5$, $c_{21} = 0.5$, $c_{22} = 1$ and homogeneous ambiguity sets, i.e, 
$\varepsilon_1 = 0.1$, $\varepsilon_2 = 0.1$, $p_1 = 1$, $p_2 = 1$. Considering heterogeneous ambiguity sets
($\varepsilon_1= 2$,
$\varepsilon_2 = 0.1$,
$p_1 = 1$,
$p_2 = -4$), we note that the equilibrium point changes. Thus, our setting generalizes standard distributionally robust optimization to games of agents with self-interested objectives that can be affected by different distributions, have their own individual data and select their own preferred risk-aversion via their Wasserstein radius. This example shows how different radii and different samples, can introduce information asymmetry in the problem at hand, thus significantly changing the set of equilibria.

\section{Wasserstein Distributionally robust  quadratic-bilinear games}
In this work, we consider a class of quadratic-bilinear games under the 2-Wasserstein distance ($p=2$).
We impose the following assumption:
\begin{assumption} \label{function_forms}
\begin{enumerate}[label=(\roman*), itemsep = 0mm, topsep = 0mm, leftmargin = 7mm]
     \item For each $i \in \mathcal{N}$, $f_i$ is convex in $x_i$ for any given $x_{-i} \in X_{-i}$;
    \item For each $i \in \mathcal{N}$, $g_i$ has the form  $g_i(x_i, x_{-i}, \xi_i)=\xi_i^\top Q_i \xi_i + P_i(x)\xi_i$;
    \item   $P_i$ is affine in $x$, i.e. $P_i(x)=A_ix+b_i$, where $A_i \in \mathbb{R}^{m \times n\SG{N}}$ and $b_i \in \mathbb{R}^m$ for all $i \in \mathcal{N}$;
    \item There exists an orthogonal matrix $L_i$ such that $Q_i=L_i^\top D_i L_i$, where $D_i$ is a diagonal positive semidefinite matrix with sorted eigenvalues. 
  \end{enumerate}
\end{assumption}

\SG{Note that $A_i$ can be written as $A_i=(A^{(1)}_i, A^{(2)}_i, \dots,  A^{(N)}_i )$, where $A^{(j)}_i \in \mathbb{R}^{m \times n}$ corresponds to the submatrix to be multiplied with the elements of $x_j$ for each $j \in \mathcal{N}$}. The structure of function $g_i$ allows for each agent to determine individually how much they wish to penalize large deviations of the uncertain parameter, represented by the quadratic term $\xi_i^\top Q_i \xi_i $. Furthermore, the bilinear term $P_i(x)\xi_i$ models the interplay between uncertainty and decisions and is important in models where the collective decision of the agents amplifies the effects of uncertainty in the cost.   Assumption (iv) allows \( Q_i \) to be represented as \( Q_i = L_i^\top D_i L_i \), where \( L_i \) is orthogonal and \( D_i \) is a diagonal positive semidefinite matrix with sorted eigenvalues. This form leverages the benefits of orthogonal transformations and simplifies the analysis of quadratic forms, while still being general enough to cover a wide range of practical scenarios.

Considering an ambiguity set per agent as in (\ref{ambiguity_set}), we then obtain the following result:
\begin{lemma} \label{Kantorovich_Duality}
Let Assumption \ref{function_forms} hold. Fix the  Wasserstein radii $\varepsilon_i$ and consider a multi-sample $\xi_{K_i}\in \mathbb{R}^{mK_i}$ for each agent $i \in \pazocal{N}$. Then, each agent's $i \in \pazocal{N}$ problem admits the following dual reformulation: 
\begin{align} \label{Robust_formulation} 
\min_{\substack{x_i \in X_i \\ \lambda_i \geq 0}}  J_i(x_i,\lambda_i, x_{-i}),  
\end{align} 
where 
\begin{align}
J_i(x_i, \lambda_i, x_{-i})=
  f_i(x_i, x_{-i})+ \lambda_i \varepsilon^2_i+ \frac{1}{K_i}\sum_{k_i=1}^{K_i}\sup_{\xi_i \in \mathbb{R}^m}[\xi_i^\top Q_i \xi_i + P_i(x)\xi_i-\lambda_i\|\xi_i-\xi_i^{(k_i)}\|^2] .  \label{Eq_ref_lemma_2_1} 
 \end{align}
\end{lemma}
\emph{Proof}: The proof follows by application of the Kantorovich duality \cite{Kantorovich_1958}. \hfill $\blacksquare$

We call the collection of the coupled optimization problems above for all agents $i \in \mathcal{N}$ as game $\bar{G}$. Note that this reformulation has an additional dual variable $\lambda_i$ for each $i \in \mathcal{N}$ that corresponds to the Lagrange multiplier associated with each individual Wasserstein constraint. Through this reformulation, a distributionally robust Nash equilibrium problem can be recast as an augmented robust Nash equilibrium problem. To connect the solutions of $G$ and $\bar{G}$, we first provide the definition of the robust Nash equilibrium (RNE) for game $\bar{G}$ as follows: 

\begin{definition} \label{Nash_definition_reformulation}
	A decision vector $(x^\ast, \lambda^\ast)$, where $\lambda^\ast= \col((\lambda_i)_{i \in \mathcal{N}})$ is a RNE of $G$ if 
	 \begin{align}
		\forall  i \in \pazocal{N}: (x^\ast_i, \lambda_i^\ast) \in \argmin_{x_i \in X_i , \lambda_i \geq 0}  J_i(x_i, \lambda_i, x^\ast_{-i}), \nonumber 
	\end{align}
 with $J_i$ as in (\ref{Eq_ref_lemma_2_1}). \hfill $\square$
\end{definition}

The following lemma establishes the connection between the set of DRNE of $G$ and the set of RNE of $\bar{G}$ defined as follows:  
\begin{lemma} \label{reformulation}
	Let $(x^\ast, \lambda^\ast)$ be an RNE of $\bar{G}$ in (\ref{Robust_formulation}). Then, $x^\ast$  is a DRNE of $G$ in (\ref{DR_formulation}). \hfill $\square$
\end{lemma}
\emph{Proof: For a given $x_{-i}^\ast \in X_{-i}$, since $(x^\ast, \lambda^\ast)$ is an RNE of $\bar{G}$, we have
\begin{align}
	& \max_{\mathbb{Q}_i \in \mathcal{P}_i} \{  f_i(x^\ast_i, x^\ast_{-i})+ \mathbb{E}_{\xi_i \sim \mathbb{Q}_i}[\xi_i^\top Q_i \xi_i + P_i(x^\ast_i, x^\ast_{-i})\xi_i \}  \nonumber \\
	&=\min_{\lambda _i \geq 0}f_i(x^\ast_i, x^\ast_{-i})+ \lambda_i \varepsilon^2_i+ \frac{1}{K_i}\sum_{k_i=1}^{K_i}\sup_{\xi_i \in \mathbb{R}^m}[\xi_i^\top Q_i \xi_i + P_i(x^\ast_i, x^\ast_{-i})\xi_i-\lambda_i\|\xi_i-\xi_i^{(k_i)}\|^2] \nonumber \\
	& = f_i(x^\ast_i, x^\ast_{-i})+ \lambda^\ast_i \varepsilon^2_i+ \frac{1}{K_i}\sum_{k_i=1}^{K_i}\sup_{\xi_i \in \mathbb{R}^m}[\xi_i^\top Q_i \xi_i + P_i(x^\ast_i, x^\ast_{-i})\xi_i-\lambda^\ast_i\|\xi_i-\xi_i^{(k_i)}\|^2] \nonumber \\
	&\leq \min_{x_i \in X_i, \lambda_i \geq 0}f_i(x_i, x^\ast_{-i})+ \lambda_i \varepsilon^2_i+ \frac{1}{K_i}\sum_{k_i=1}^{K_i}\sup_{\xi_i \in \mathbb{R}^m}[\xi_i^\top Q_i \xi_i + P_i(x_i, x^\ast_{-i})\xi_i-\lambda_i\|\xi_i-\xi_i^{(k_i)}\|^2] \nonumber \\
	&=   \min_{x_i \in X_i } \max_{\mathbb{Q}_i \in \mathcal{P}_i} \{  f_i(x_i, x^\ast_{-i})+ \mathbb{E}_{\xi_i \sim \mathbb{Q}_i}[\xi_i^\top Q_i \xi_i + P_i(x_i, x^\ast_{-i})\xi_i] \},
\end{align}
where  the  inequality holds from Definition \ref{Nash_definition_reformulation}.  \hfill $\blacksquare$ \\}   

Note that the inverse direction does not necessarily hold, as one should determine an appropriate value for $\lambda^\ast$. From Lemma \ref{reformulation} we can instead solve game $\bar{G}$ and obtain the solution $(x^\ast, \lambda^\ast)$ and, from this solution, select $x^\ast$ as the DRNE of our original problem $G$. \GP{To achieve this, we impose the following standing assumption:}
\begin{assumption} \label{non_empty}
	The set of RNE of $\bar{G}$ in (\ref{Robust_formulation}) is non-empty.  \hfill $\square$
\end{assumption} 
The non-emptiness of the set of RNE of $\bar{G}$ then directly implies the non-emptiness of the set of DRNE of game $G$. To solve the inner maximization over the uncertain parameter $\xi_i$, we show  that the class of games that satisfies Assumption  \ref{function_forms} can be exploited to provide a  finite-dimensional formulation, without the use of an epigraphic reformulation. The following theorem leverages the structure of the problem to obtain a more computationally efficient reformulation thus circumventing those challenges.   \par 

\begin{thrm} \label{theorem1}
	Under Assumption \ref{function_forms},  $G$ admits the reformulation
	\begin{align}
		G_R: \forall i \in \mathcal{N}:	\min_{x_i \in X_i, \lambda_i > \lambda_{max}(Q_i)}&\{f_i(x_i, x_{-i})+ \lambda_i\left(\varepsilon^2_i-\frac{1}{K_i}\sum_{k_i=1}^{K_i} (\xi^{(k_i)}_i)^\top\xi^{(k_i)}_i\right)+ \nonumber \\
		&+\frac{1}{4K_i}\sum_{k_i=1}^{K_i}\tilde{W}^{(k_i)}(x_i,x_{-i}, \lambda_i)^\top \tilde{Q}_i(\lambda_i) \tilde{W}^{(k_i)}(x_i,x_{-i}, \lambda_i)\}, \nonumber
	\end{align}
	where $\tilde{Q}_i(\lambda_i)=\text{diag}(\frac{1}{\lambda_i-\lambda_{max}(D_i)}, \dots, \frac{1}{\lambda_i-\lambda_{min}(D_i)})$ and $\tilde{W}^{(k_i)}(x_i,x_{-i}, \lambda_i)=\tilde{P}_i(x_i, x_{-i})+2\lambda_i\ \tilde{\xi}^{(k_i)}_i$, $\tilde{P}_i(x_i, x_{-i})=L_iP_i(x_i, x_{-i})$ and $\tilde{\xi}^{(k_i)}_i=L_i\xi^{(k_i)}_i$. \hfill $\square$ 
\end{thrm}
\emph{Proof}:
For each agent $i \in \mathcal{N}$ it holds that:
\begin{align}
	& \frac{1}{K_i}\sum_{k_i=1}^{K_i}\sup_{\xi_i \in \mathbb{R}^m}[\xi_i^\top Q_i \xi_i + P_i(x_i, x_{-i})\xi_i-\lambda_i\|\xi_i-\xi_i^{(k_i)}\|^2] \nonumber \\
	& = \frac{1}{K_i}\sum_{k_i=1}^{K_i}\sup_{\xi_i \in \mathbb{R}^m}[\xi_i^\top Q_i \xi_i + P_i(x_i, x_{-i})\xi_i-\lambda_i(\xi_i^\top \xi_i -2 \xi_i^\top \xi^{(k_i)}_i+ (\xi^{(k_i)}_i)^\top \xi^{(k_i)}_i)]  \nonumber \\
	& = \frac{1}{K_i}\sum_{k_i=1}^{K_i} \left(-\lambda_i (\xi^{(k_i)}_i)^\top \xi^{(k_i)}_i+\sup_{\xi_i \in \mathbb{R}^m}[\xi_i^\top (Q_i-\lambda_i I_{m}) \xi_i + (P_i(x_i, x_{-i})+ 2 \lambda_i  \xi^{(k_i)}_i)^\top \xi_i]  \right) \nonumber 
\end{align}

Since for each $i \in \mathcal{N}$, $Q_i$ is diagonalizable, there exists matrix $L_i$ such that  $Q_i=L_i^\top D_i L_i$,  where $D_i$ is a diagonal matrix, whose eigenvalues decrease along the diagonal. Denote the maximum eigenvalue of $D_i$ by $\lambda_{max}(D_i)$ and the minimum eigenvalue of $D_i$ by $\lambda_{min}(D_i)$. As such,  the following equalities hold:

\begin{align}
	& \frac{1}{K_i}\sum_{k_i=1}^{K_i} \left(-\lambda_i (\xi^{(k_i)}_i)^\top \xi^{(k_i)}_i+\sup_{\xi_i \in \mathbb{R}^m}[\xi_i^\top (Q_i-\lambda_i I_{m}) \xi_i + (P_i(x_i, x_{-i})+2 \lambda_i  \xi^{(k_i)}_i)^\top \xi_i]  \right) \nonumber \\
	&= \frac{1}{K_i}\sum_{k_i=1}^{K_i} \left(-\lambda_i (\xi^{(k_i)}_i)^\top \xi^{(k_i)}_i+\sup_{\xi_i \in \mathbb{R}^m}[\xi_i^\top (D_i-\lambda_i I_{m}) \xi_i +[L_i P_i(x_i, x_{-i})+2 \lambda_i  L_i\xi^{(k_i)}_i]^\top \xi_i]  \right)  \nonumber  \\
   &= \frac{1}{K_i}\sum_{k_i=1}^{K_i} \left(-\lambda_i (\xi^{(k_i)}_i)^\top \xi^{(k_i)}_i+\sup_{\xi_i \in \mathbb{R}^m}[\xi_i^\top (D_i-\lambda_i I_{m}) \xi_i  +\tilde{W}^{(k_i)}(x_i,x_{-i}, \lambda_i)^\top \xi_i]  \right)  \label{intermediate_der}
\end{align}
Consider now $G^{(k_i)}_i(x_i, x_{-i}, \lambda_i,\xi_i )=\xi_i^\top (D_i-\lambda_i I_{m}) \xi_i +[\tilde{P}_i(x_i, x_{-i})\xi_i+2 \lambda_i  \tilde{\xi}^{(k_i)}_i)]^\top \xi_i$. Due to the presence of the supremum in (\ref{intermediate_der}), we wish to study for which value of the uncertainty $\xi_i$ we achieve the maximum value for $G^{(k_i)}_i(x_i, x_{-i}, \lambda_i,\xi_i )$. This maximum value will be parametrized by the corresponding sample $\xi^{(k_i)}_i$. We distinguish between two different cases:
\begin{enumerate}[label=(\roman*)]
	\item   For $\lambda_i > \lambda_{max}(Q_i)$ we note that $(D_i-\lambda_i I_m)^{-1}$ is negative semidefinite. Thus, given other agents' decisions $x_{-i}$, the resulting cost function is concave which yields the solution
 \begin{align}
 \sup_{\xi_i \in \mathbb{R}^m}\xi_i^\top (D_i-\lambda_i I_{m}) \xi_i +[\tilde{P}_i(x_i, x_{-i})+2 \lambda_i  \tilde{\xi}^{(k_i)}_i)]^\top \xi_i=G^{(k_i)}_i(x_i, x_{-i}, \lambda_i,\xi^\ast_i ),
 \end{align}  where $\xi^\ast_i$ is obtained by the first order optimality condition $\nabla_{\xi_i^\ast}G^{(k_i)}_i(x_i, x_{-i}, \lambda_i,\xi_i )=0$. As such, the maximum is attained at $\xi^\ast_i=\frac{1}{2}(\lambda_i I_m-D_i)^{-1}(\tilde{P}_i(x_i, x_{-i})+2\lambda_i\tilde{\xi}_i^{(k_i)})$ with optimal value:
	\begin{align}
		&G^{(k_i)}_i(x_i, x_{-i}, \lambda_i,\xi^\ast_i )=(\xi^\ast_i)^\top (D_i-\lambda_i I_{m}) \xi^\ast_i +[\tilde{P}_i(x_i, x_{-i})+2 \lambda_i  \tilde{\xi}^{(k_i)}_i)]^\top \xi_i^\ast \nonumber \\
		&=  \frac{1}{4}(\tilde{P}_i(x_i, x_{-i})+ 2\lambda_i \tilde{\xi}_i^{(k_i)})^\top (\lambda_i I_m-D_i)^{-1}(\tilde{P}_i(x_i, x_{-i})+2\lambda_i \tilde{\xi}^{(k_i)}))
	\end{align}
	\item  For $ \lambda_i \in [0, \lambda_{max}(Q_i))$, we note that $D_i-\lambda_i I_m$ is positive semidefinite, hence the cost function of the inner maximization problem is convex in $\xi_i$, which implies that $\sup_{\xi_i \in \mathbb{R}^m}\xi_i^\top (D_i-\lambda_i I_{m}) \xi_i +[\tilde{P}_i(x_i, x_{-i})+2 \lambda_i  \tilde{\xi}^{(k_i)}_i)]^\top \xi_i= \infty$.
\end{enumerate}
As such, given the agents' decisions $x_{-i}$ each agent $i \in \mathcal{N}$ solves 
\begin{align}
		\min_{x_i \in X_i, \lambda_i > \lambda_{max}(Q_i)}&\{f_i(x_i, x_{-i})+ \lambda_i\left(\varepsilon^2_i-\frac{1}{K_i}\sum_{k_i=1}^{K_i} (\xi^{(k_i)}_i)^\top\xi^{(k_i)}_i\right)+ \nonumber \\
		&+\frac{1}{4K_i}\sum_{k_i=1}^{K_i}\tilde{W}^{(k_i)}(x_i,x_{-i}, \lambda_i)^\top \tilde{Q}_i(\lambda_i) \tilde{W}^{(k_i)}(x_i,x_{-i}, \lambda_i)\}, \nonumber
	\end{align}
 where $\tilde{Q}_i(\lambda_i) =(\lambda_i I_m-D_i)^{-1}$.
Then, the connection between the games $G$ and $\bar{G}$ as established in Lemma \ref{Kantorovich_Duality} and their corresponding solutions in Lemma \ref{reformulation} concludes the proof. \hfill $\blacksquare$ \\

\subsection{Reformulation as a Data-driven Variational Inequality Problem}
 \SG{In this section, we establish the connection of the Nash equilibria of $G_R$ with the solutions of a variational inequality (VI) problem. For the ease of the reader, we define the notion of a Nash equilibrium for a general game. } 
 
\begin{definition} \label{NE_defin}
Consider the following game:
 \begin{align}
\forall \ i \in \mathcal{N}: \min\limits_{z_i \in Z_i}J_i(z_i, z_{-i}), \label{game2}
\end{align}
A point $z^\ast=(z_i^\ast, z_{-i}^\ast) \in Z=\prod_{i=1}^N Z_i$ is Nash equilibrium (NE) of (\ref{game2}) if, given $x^\ast_{-i}$, the following condition holds:
\begin{align}
J_i(z^\ast_i, z^\ast_{-i}) \leq J_i(z_i, z^\ast_{-i}), \nonumber
\end{align}
for all $z_i \in Z_i$ and for all $i \in \mathcal{N}$. \hfill $\square$  
\end{definition}
The following statement then holds:

 \begin{proposition} \label{game_VI_equivalence}
 Consider the following game:
 \begin{align}
\forall \ i \in \mathcal{N}: \min\limits_{z_i \in Z_i}J_i(z_i, z_{-i}), \label{game}
\end{align}
where $J_i$ is convex on $Z_i$ for any $z_{-i} \in Z_{-i}$ and $Z=\prod_{i=1}^N Z_i$ is convex and closed. Furthermore, consider the following variational inequality problem:
	\begin{align}
		F^\top(z^\ast)(z-z^\ast) \geq 0,  \ \forall \ z \in Z \cap V(z^\ast) \text{ for all }  i \in \mathcal{N}.  \nonumber
	\end{align}
	where 
	$F(z)=\col((F_i(z))_{i \in \mathcal{N}})$ with $F_i(z)=\nabla_{z_i}J_i(z_i, z_{-i})$ is a (possibly nonmonotone) mapping and $V(z^\ast)$ is a small enough convex neighbourhood around $z^\ast$.
	Then, any local solution $z^\ast$ of the VI is a Nash equilibrium of (\ref{game}).
\end{proposition}
\emph{Proof}:  This result is a direct extension of the proofline of Proposition 1.4.2 in \cite{Pang1} for a nonmonotone mapping defined over a small enough convex neighbourhood $V(z^\ast)$ of the solution. \hfill $\blacksquare$

\SG{Returning to game $G_R$, note that, according to Definition \ref{NE_defin},  a point   $(x^\ast, \lambda^\ast) \in X \times \prod_{i=1}^N(\lambda_{max}(Q_i)), \infty)$ is a Nash equilibrium of $G_R$ if, given $x^\ast_{-i} \in X_{-i}$, the following condition holds: 
\begin{align}
  &f_i(x^\ast)+ \lambda^\ast_i\left(\varepsilon^2_i-\frac{1}{K_i}\sum_{k_i=1}^{K_i} (\xi^{(k_i)}_i)^\top\xi^{(k_i)}_i\right)+\frac{1}{4K_i}\sum_{k_i=1}^{K_i}\tilde{W}^{(k_i)}(x^\ast, \lambda_i^\ast)^\top \tilde{Q}_i(\lambda_i^\ast) \tilde{W}^{(k_i)}(x^\ast, \lambda_i^\ast) \leq \nonumber  \\
&f_i(x_i, x^\ast_{-i})+ \lambda_i\left(\varepsilon^2_i-\frac{1}{K_i}\sum_{k_i=1}^{K_i} (\xi^{(k_i)}_i)^\top\xi^{(k_i)}_i\right)\!\!+\frac{1}{4K_i}\sum_{k_i=1}^{K_i}\tilde{W}^{(k_i)}(x_i,x^\ast_{-i}, \lambda_i)^\top \tilde{Q}_i(\lambda_i) \tilde{W}^{(k_i)}(x_i,x^\ast_{-i}, \lambda_i) \nonumber
  \end{align}
 for all $(x_i, \lambda_i) \in X_i \times (\lambda_{max}(Q_i), \infty)$ \label{equilibrium_condition}}. \par 
\SG{Finally, from the proofline of Theorem \ref{theorem1}, it immediately follows that the set of NE of $G_R$ coincides with the set of RNE of (\ref{Robust_formulation})}. Let us now denote by $z_i=(x_i, \lambda_i) \in \mathbb{R}^{(n+1)N}$ the collection of the decision vector and the Lagrange multiplier for all $i \in \mathcal{N}$ and $z=\col((z_i)_{i \in \mathcal{N}})$. Furthermore, denote the feasible set  $Z=\{ z \in  \mathbb{R}^{(n+1)N}: x_i \in X_i, \lambda_i \geq \lambda_{max}(Q_i)+\zeta_i, \ \forall i \in \mathcal{N} \}$, where $\zeta_i$ is an arbitrarily small positive parameter, ensuring that the local constraint set is closed and thus game $G_R$, where $\lambda_i> \lambda_{\max}(Q_i)$, can be solved with any a priori defined accuracy. An exact solution is obtained when $\zeta_i \rightarrow 0^+$.
The following lemma then holds:
\begin{lemma}
	\SG{A solution  $z^\ast$ of the VI problem with mapping $F(z)=\col(F_i(z))_{i \in \mathcal{N}}$, where
 \begin{align} \label{reformulation_VI}
		F_i(z)=
		\begin{pmatrix}\nabla_{x_i}f_i(x_i, x_{-i})+\dfrac{1}{2K_i}\sum\limits_{k_i=1}^{K_i} (A^{(i)}_i)^\top\tilde{Q}_i(\lambda_i)\tilde{W}^{(k_i)}(x_i,x_{-i}, \lambda_i)  \\
			\varepsilon_i^2  -\dfrac{1}{K_i}\sum\limits_{k_i=1}^{K_i} \|\xi_i^{(k_i)}\|^2+\dfrac{1}{4K_i}\sum\limits_{k_i=1}^{K_i}4\xi^{(k_i) \top}_{i} 	\tilde{Q}_i(\lambda_i)\tilde{W}^{(k_i)}(x, \lambda_i)- \| \tilde{W}^{(k_i)}(x, \lambda_i)\|^2_{\frac{d \tilde{Q}_i}{d \lambda_i}}  \\
		\end{pmatrix}. 
	\end{align} 
  over $Z \cap V(z^\ast)$, with  $V(z^\ast)$ being a convex local neighbourhood around $z^\ast$, is a Nash equilibrum  of $G_R$ over $Z$}.

\end{lemma}

\emph{Proof}:  The proof \GP{follows by} adaptation of Proposition \ref{game_VI_equivalence} by taking the pseudogradient of game $G_R$ considering that $\tilde{Q}_i$ is a diagonal matrix, hence $\frac{d\tilde{Q}_i(\lambda_i)}{d \lambda_i}$ is obtained by differentiating the corresponding diagonal elements. \hfill $\blacksquare$
\
\subsection{\RB{Challenges of distributionally robust heterogeneous  games compared to DRO}}

Let us illustrate the challenges that our game-theoretic problem entails compared to the standard DRO. 
Let us first define the associated distributionally robust optimization problem as follows:

\begin{equation}
\min_{x \in X} \max_{\mathbb{Q} \in \mathcal{P}_{\epsilon}(\hat{\mathbb{P}}_K)} \mathbb{E}_{\xi \sim \mathbb{Q}}[x^\top Qx+\xi^\top Q \xi + P(x)^\top\xi], \label{eq:DRO}
\end{equation}

First, note that the optimization problem in (\ref{eq:DRO}) is only a special case of the original game considered in (\ref{DR_formulation}). In particular, the setting in (\ref{DR_formulation}) considers more than one players, who make self-interested decisions based on  their own \emph{private samples} and \emph{different} Wasserstein radii. 
The interdependence of each agent's cost to the decisions of the other agents significantly affects the resulting decision, unlike optimization where this decision is collective and only nature is viewed as an adversarial player. \GP{This implies that the pseudogradient mapping might be non-monotone; On the contrary, in a standard optimization problem, the gradient map is always monotone.}

\GP{For example, consider the following distributionally robust game with 2 agents $i \in \{1, 2\}$ taking scalar decisions. Specifically, agent 1 aims at solving:}
\begin{equation}
  \min_{x_1 \in X_1} \max_{\mathbb{Q}_1 \in \mathcal{P}_1} \mathbb{E}_{\xi_1 \sim \mathbb{Q}_1}[q_{11}x^2_1+ q_{12}x_1 x_2+ Q_1 \xi^2_1 + (a_1x_1+a_2x_2)\xi_1] , \nonumber
\end{equation}
while agent 2 aims at solving:
\begin{align}
	\min_{x_2 \in X_2} \max_{\mathbb{Q}_2 \in \mathcal{P}_2} \mathbb{E}_{\xi_2\sim \mathbb{Q}_2}[q_{22}x^2_2+ q_{21}x_1 x_2+ Q_2 \xi^2_2 + (b_1x_1+b_2x_2)\xi_2] , \nonumber
\end{align}

By Theorem 1, the problem above can be written as:
\begin{equation} \label{eq:game_2d}
\forall \ i \in \mathcal{N}:	
		\min_{(x_i, \lambda_i) \in X_i \times \mathbb{R}^m_{\geq 0}} J_i(x_i,x_{j}, \lambda_i),
	\end{equation}
where
\begin{equation}
J_1(x_1,x_2, \lambda_1)=q_{11}x^2_1+ q_{12}x_1 x_2+ \lambda_1\left(\epsilon^2_1-\frac{1}{K_1}\sum_{k_1=1}^{K_1} (\xi^{(k_1)}_1)^2\right)+\frac{1}{4K_1}\sum_{k_1=1}^{K_1} \frac{(a_1x_1+a_2x_2+2 \lambda_1 \xi^{(k_1)}_1)^2}{\lambda_1-\lambda_{max,1}}, \nonumber 
\end{equation}
and
\begin{equation}
J_2(x_1,x_2, \lambda_2)=q_{22}x^2_2+ q_{21}x_1 x_2+ \lambda_2\left(\epsilon^2_2-\frac{1}{K_2}\sum_{k_2=1}^{K_2} (\xi^{(k_2)}_2)^2\right)+\frac{1}{4K_2}\sum_{k_2=1}^{K_2} \frac{(b_1x_1+b_2x_2+2 \lambda_2 \xi^{(k_2)}_2)^2}{\lambda_2-\lambda_{max,2}}. \nonumber
\end{equation}
To study the monotonicity of the VI mapping associated with (\ref{eq:game_2d}), we can directly study whether the game Jacobian is positive semidefinite. The following lemma shows that for this class of games, even if $K_1=K_2$ (note that only the number is the same, thus still allowing for different samples per agent), the Jacobian is more complicated than the calculated one for standard optimization. 
\begin{lemma}  \label{lem:Jacobian}
Consider the functions \begin{align}
	&G^{(k_1)}_1(x_1,x_2, \lambda_1)=2\xi^{(k_1)}_1\lambda_{1, \max}+a_1 x_1+a_2 x_2, \nonumber  \\
	&G^{(k_2)}_2(x_1,x_2, \lambda_2)=2\xi^{(k_2)}_2\lambda_{2, \max}+b_1x_1+b_2x_2, \nonumber \\
    &\Delta_1=\lambda_1-\lambda_{1, \max} \nonumber \\
    &\Delta_2=\lambda_2-\lambda_{2, \max}, \label{2_player_game}
\end{align}
corresponding to each sample $k_1 \in \{1, \dots, K_1\}$ and $k_2 \in \{1, \dots, K_2\}$. Furthermore, consider for simplicity that $K_1=K_2=K$ and $q_{11}=q_{12}=q_{22}=q_{21}=0$. The symmetrized Jacobian $J$ is then given by
\begin{equation}
J=\frac{1}{4K}\sum_{k=1}^K J^{(k)}_s, \nonumber
\end{equation}
 where
\begin{equation}
J^{(k)}_s = \begin{bmatrix}
	\frac{2a_1^2}{\Delta_1} & \frac{2a_1G^{(k)}_1(x, \lambda_1)}{\Delta_1^2} & \frac{a_1a_2}{\Delta_1} + \frac{b_1b_2}{\Delta_2} & 	\frac{b_1G^{(k)}_2(x, \lambda_2)}{\Delta_2^2} \\
	\frac{2a_1G^{(k)}_1(x, \lambda_1)}{\Delta_1^2} & \frac{2G^{(k)}_1(x, \lambda_1)^2}{\Delta_1^3} & \frac{a_2G^{(k)}_1(x, \lambda_1)}{\Delta_1^2} & 0 \\
	\frac{a_1a_2}{\Delta_1} + \frac{b_1b_2}{\Delta_2} &\frac{a_2G^{(k)}_1(x, \lambda_1)}{\Delta_1^2} & \frac{2b_2^2}{\Delta_2} & \frac{2b_2G^{(k)}_2(x, \lambda_2)}{\Delta_2^2} \\
	\frac{b_1G^{(k)}_2(x, \lambda_2)}{\Delta_2^2} & 0 & \frac{2b_2G^{(k)}_2(x, \lambda_2)}{\Delta_2^2} & \frac{2G^{(k)}_2(x, \lambda_2)^2}{\Delta_2^3} \label{Jacobian_game}
\end{bmatrix}.
\end{equation} 
\end{lemma}
\emph{Proof}: The proof is obtained by performing standard algebraic operations on the defined Jacobian of the corresponding VI, considering the functions $G^{(k_1)}$ and  $G^{(k_2)}$. We then use the assumption $K_1=K_2=K$ to pull a common sum outside and group the samples of the two agents in exhaustive pairs. \hfill $\blacksquare$

Based on  Lemma \ref{lem:Jacobian}, an important observation can then be made: When we consider \emph{only} one agent, we derive as a special case of our original problem a distributionally robust optimization program with Jacobian:
\begin{equation}
J^{(k)}_s = \begin{bmatrix}
	\frac{2a_1^2}{\Delta_1} & \frac{2a_1G^{(k_1)}_1(x, \lambda_1)}{\Delta_1^2} \\
	\frac{2a_1G^{(k_1)}_1(x, \lambda_1)}{\Delta_1^2}  &
    \frac{2G^{(k_1)}_1(x, \lambda_1)^2}{\Delta_1^3} 
\end{bmatrix},
\end{equation}

For any $k$, $J_s^{(k)}$ can then be shown to be positive semi-definite. Thus $J$ is positive semi-definite as the summation of positive semi-definite matrices. In contrast,  the Jacobian (\ref{Jacobian_game}) for the game in (\ref{eq:DRO}) might not be positive semi-definite which implies non-monotonicity. For example, the game in (\ref{2_player_game}) with fixed decisions $x_1=2$, $x_2=2$, $\lambda_1=\lambda_{1,\max}+2$ and  $\lambda_2=\lambda_{2,\max}+2$, and parameters $\lambda_{1,\max} = 0.5$, $\lambda_{2,\max} = 0.5$, $a_1 = 1$, $a_2 = 0.5$,  
$b_1 = 1$, $b_2 = 0.5$, has Jacobian with at least one negative eigenvalue meaning that the corresponding pseudogradient is not monotone.

Even though the mapping can in general be nonmonotone, we illustrate how equilibria can still be computed efficiently based both on the structure of our problem, as provided by Theorem \ref{theorem1}, and the equilibrium seeking algorithms we propose in the next section. Specifically, leveraging the structure of our reformulation, we avoid introducing epigraphic variables that render both the decision and the constraints dependent on the size of the data \cite{netessine_wasserstein_2019}, \cite{mohajerin_esfahani_data-driven_2018}.  Thus, through Theorem \ref{theorem1}, we can obtain  reformulations for the class of heterogeneous data-driven Wasserstein distributionally robust games in (\ref{DR_formulation}) that scale better with data compared to the use of epigraph forms. In the next section, we assess the computational performance of our theoretical results through an illustrative example and a risk-aware portfolio allocation game, which takes into account behavioural coupling of the investors' decisions.

\section{Numerical Simulations}
\begin{algorithm}[t]
\caption{Hybrid DRNE seeking algorithm ({Hybrid-Alg})\cite{Reza_2024}}
\label{{Hybrid-Alg}2}
\begin{algorithmic}[1]
\REQUIRE Choose $x^0$, $x^1$, $\tau_0 > 0$, $\bar\tau \gg 0$, $\alpha = (1,\frac{1+\sqrt{5}}{2}]$, $\theta_0 = 1$, $\rho = \dfrac{1}{\alpha} + \dfrac{1}{\alpha^2}$, $\bar \phi \gg \frac{1+\sqrt{5}}{2}$, $\text{sum}_0 ^1 = 0$, $\text{sum}_0 ^2 = 0$, flg = 1.
\STATE \textbf{For} {$k = 0,1,2,\ldots$} \textbf{do}
\STATE Find the stepsize:
\qquad \qquad \qquad \[\tau_k = \min\left\{\rho \tau_{k-1}, \dfrac{\alpha\theta_{k-1}}{4\tau_{k-1}} \dfrac{\|x^k - x^{k-1}\|^2}{\|F(x^k) - F(x^{k-1})\|^2}, \Bar{\tau}\right\}\]
\STATE $\bar{x}^{k} = \dfrac{(\phi_k - 1)x^k + \bar{x}^{k-1}}{\phi_k}$
\STATE Update the next iteration:\\
\qquad \qquad \qquad $x^{k+1} = \text{prox}_{\tau_k g}(\bar{x}^k - \tau_k F(x^k))$
\STATE Update: \qquad $\theta_{k+1} = \dfrac{\alpha\tau_k}{ \tau_{k-1}}$
\STATE compute the following summations with $\phi_{k+1} = \bar \phi$:\\
\qquad \qquad \qquad $\text{sum}_{k+1}^1 = \text{sum}_k^1$ + \eqref{eq1:hyb-alg}\\
\qquad \qquad \qquad $\text{sum}_{k+1}^2 = \text{sum}_k^2$ + \eqref{eq2:hyb-alg}
\IF{ ($\text{sum}_{k+1}^1 \leq 0 \,\,\, \land \,\,\, \text{flg} = 1$) \,\,$\lor$\,\, ($\text{sum}_{k+1}^2 \leq 0 \,\,\, \land \,\,\, \text{flg} = 0$)}
    \STATE $\phi_{k+1} = \bar\phi$, $\text{flg} = 1$ 
\ELSE
\IF{$\text{flg} = 1$} 
    \STATE $x^{k+1} = x^k$, $\bar x^k = \bar x^{k-1}$ 
    \STATE $\phi_{k+1} = \alpha$, $\theta_k = \theta_{k-1}$, $\tau_{k} = \tau_{k-1}$ 
    \STATE $\text{sum}_{k+1}^1 = 0$, $\text{sum}_{k+1}^2 = 0$, $\text{flg} = 0$
    \ELSE  
    \STATE $\phi_{k+1} = \alpha$ 
    \STATE $\text{sum}_{k+1}^2 = \text{sum}_k^2$ + (\eqref{eq2:hyb-alg} with $\phi_{k+1}  = \alpha$) 
    \STATE $\text{sum}_{k+1}^1 = 0$ 
\ENDIF 
\ENDIF
\end{algorithmic}
\end{algorithm}

In the simulation results, we use three algorithms to solve the variational inequality problems: 
\begin{enumerate}[label=(\roman*)]
    \item Forward-backward splitting (FB):
    \begin{align*}
    x^{k+1} = \text{prox}_{\tau g}(x^{k} - \tau F(x^k)),
    \end{align*}
    where $\tau$ is the stepsize. This method is known as a projected gradient descent algorithm as well which is the most used method in machine learning and control applications. The convergence of this method is guaranteed for strongly monotone (with a strongly monotone constant $\mu$) and Lipschitz (with a Lipschitz constant $L$) operator with $\tau \in \left(0, 2\mu/L^2\right)$ \cite[Theorem~12.4.6]{Pang1}.
    \item Adaptive golden ration algorithm (aGRAAL) \cite{malitsky_golden_2020}: 
     \begin{align*}
        \phi &\in (1,\frac{1+\sqrt{5}}{2}],\quad \rho = \dfrac{1}{\phi} + \dfrac{1}{\phi^2}\\
        \tau_k &= \min\left\{\rho \tau_{k-1}, \dfrac{\phi\theta_{k-1}}{4\tau_{k-1}} \dfrac{\|x^k - x^{k-1}\|^2}{\|F(x^k) - F(x^{k-1})\|^2}, \bar{\tau}\right\}\\
        \bar{x}^{k} &= \dfrac{(\phi - 1)x^k + \bar{x}^{k-1}}{\phi}\\
         x^{k+1} &= \text{prox}_{\tau_k g}(\bar{x}^k - \tau_k F(x^k)), \quad \theta_{k+1} = \dfrac{\phi\tau_k}{ \tau_{k-1}}
    \end{align*}
The convergence of this method is guaranteed for Lipschitz and monotone operator (with a Lipschitz constant $L$) \cite[Lemma 3]{malitsky_golden_2020}.
\item Hybrid method ({Hybrid-Alg}): \GP{The full steps of this method can be found in Algorithm \ref{{Hybrid-Alg}2}. This algorithm is presented in the recent work \cite{Reza_2024} and includes some of the authors of the present paper. However, we note that the convergence rate proof only holds for monotone mappings and is applied on a significantly less complicated case study. The application of this algorithm on our problem shows that the equilibria of even a generally nonmonotone mapping could be numerically obtained relatively fast, even faster than the golden ratio method in practice.} The Hybrid method, though similar to aGRAAL, differs significantly in the choice of the momentum parameter, as we explain later on. 
\end{enumerate}

\begin{figure}[h]
    \centering
    \subfloat[$\varepsilon = 10^{-6}$.]{\label{fig-conv:1}\includegraphics[scale=0.35]{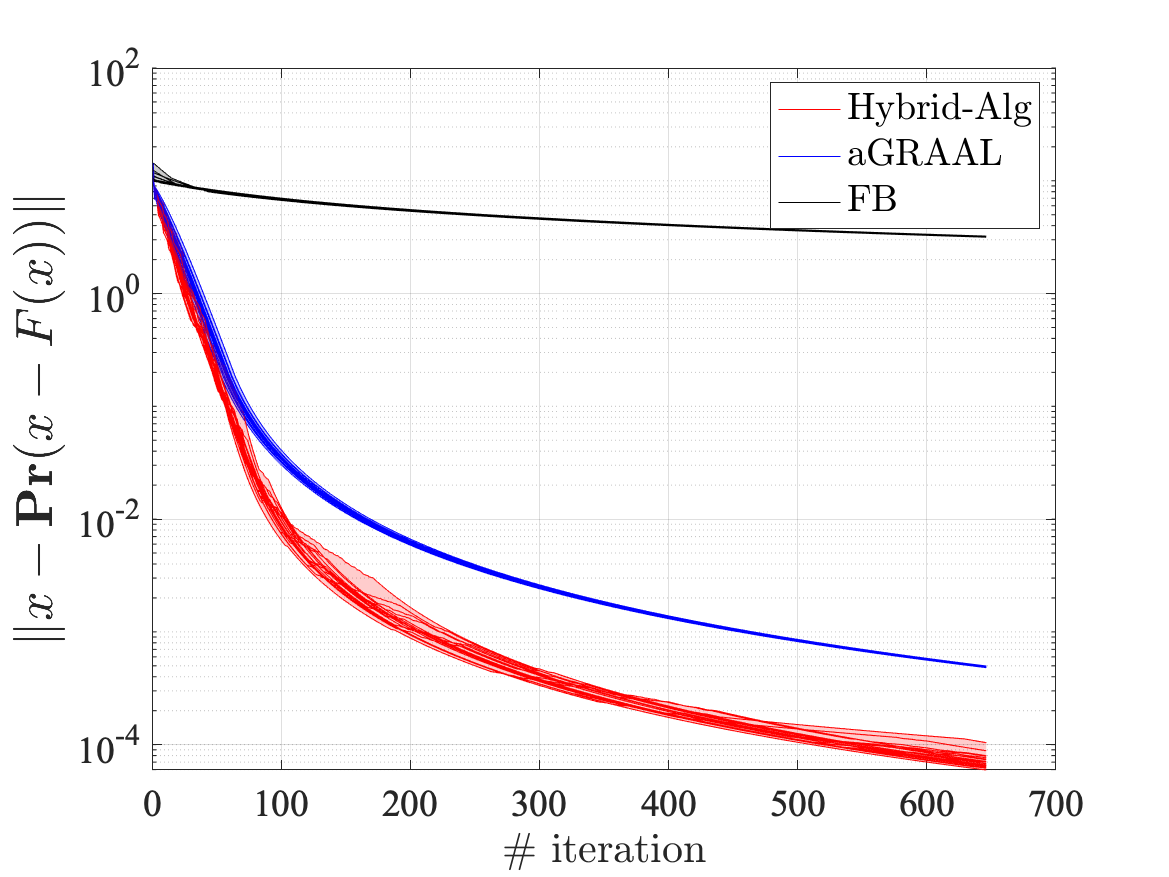}}
    \hspace{2mm}
    \subfloat[$\varepsilon = 10^{-3}$.]{\label{fig-conv:2}\includegraphics[scale=0.35]{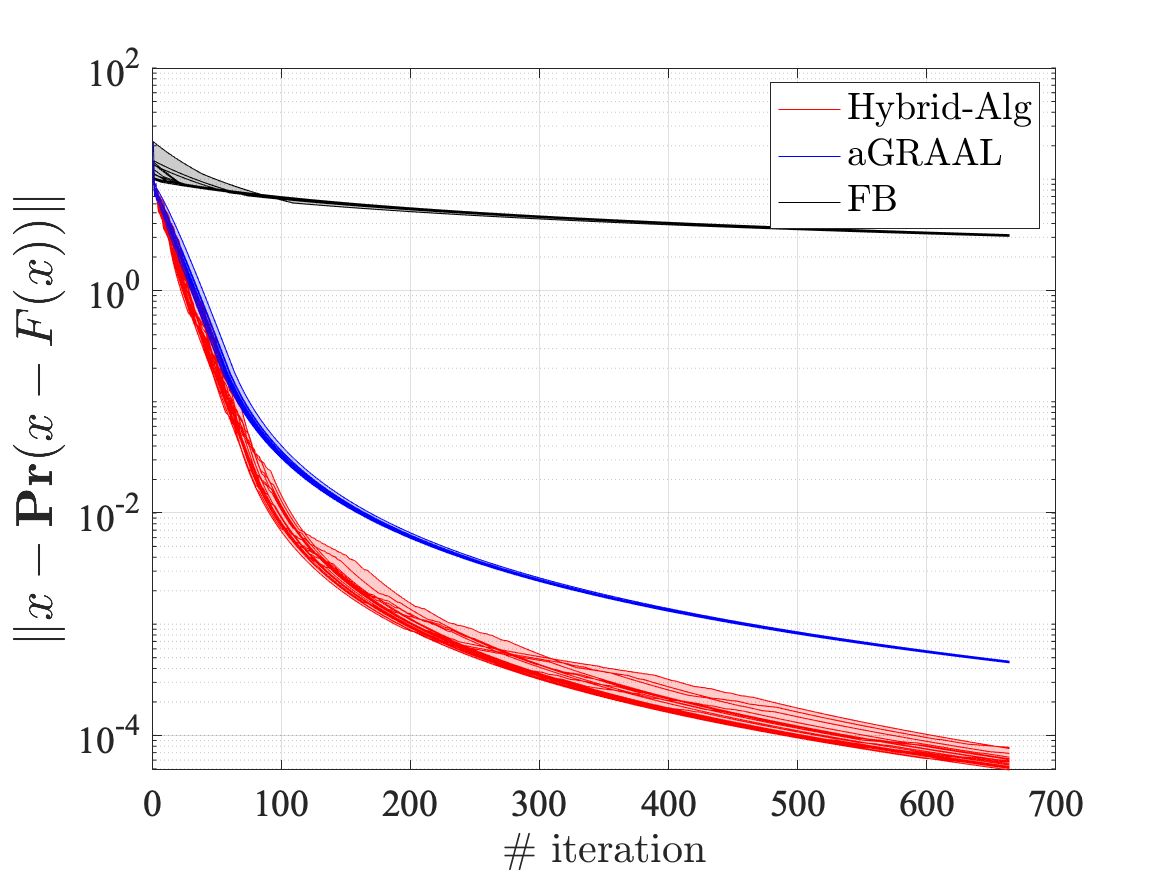}}
        \hspace{2mm}
    \subfloat[$\varepsilon = 1$.]{\label{fig-conv:5}\includegraphics[scale=0.35]{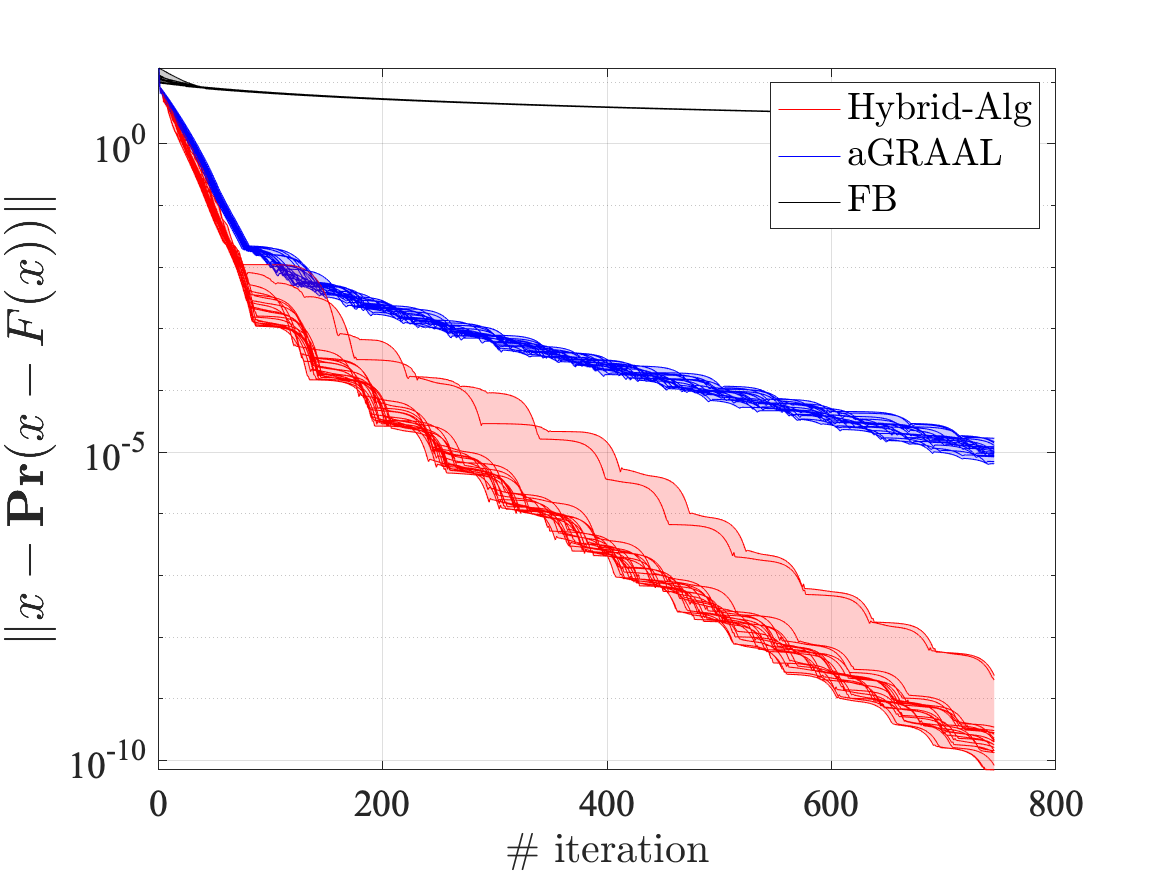}}
    \caption{Residual of variational inequality problem with different radii $\varepsilon_i$ per agent chosen according to the distribution $\varepsilon \cdot U[1, 5]$, where $\varepsilon$ takes values in $\{10^{-6}, 10^{-3}, 1\}$.}
    \label{fig-conv-radii}
\end{figure}
\begin{figure}[h] 
    \centering
    \captionsetup{justification=centering}
    \subfloat[$K_i = \text{randi}\left((10, 20),1,n\right)$.]{\label{fig-conv-smp:1}\includegraphics[scale=0.35]{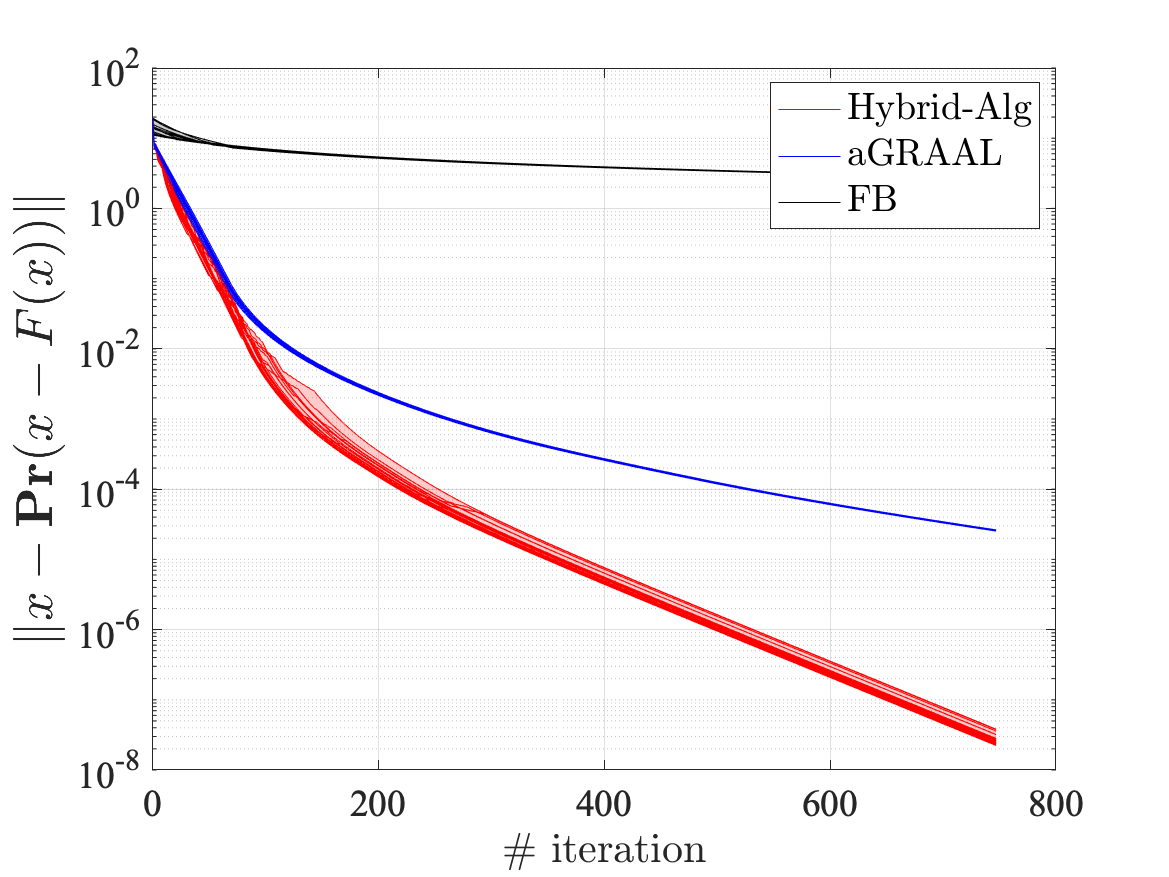}}
    \hspace{2mm}
    \subfloat[$K_i = \text{randi}\left((40, 60),1,n\right)$.]{\label{fig-conv-smp:2}\includegraphics[scale=0.35]{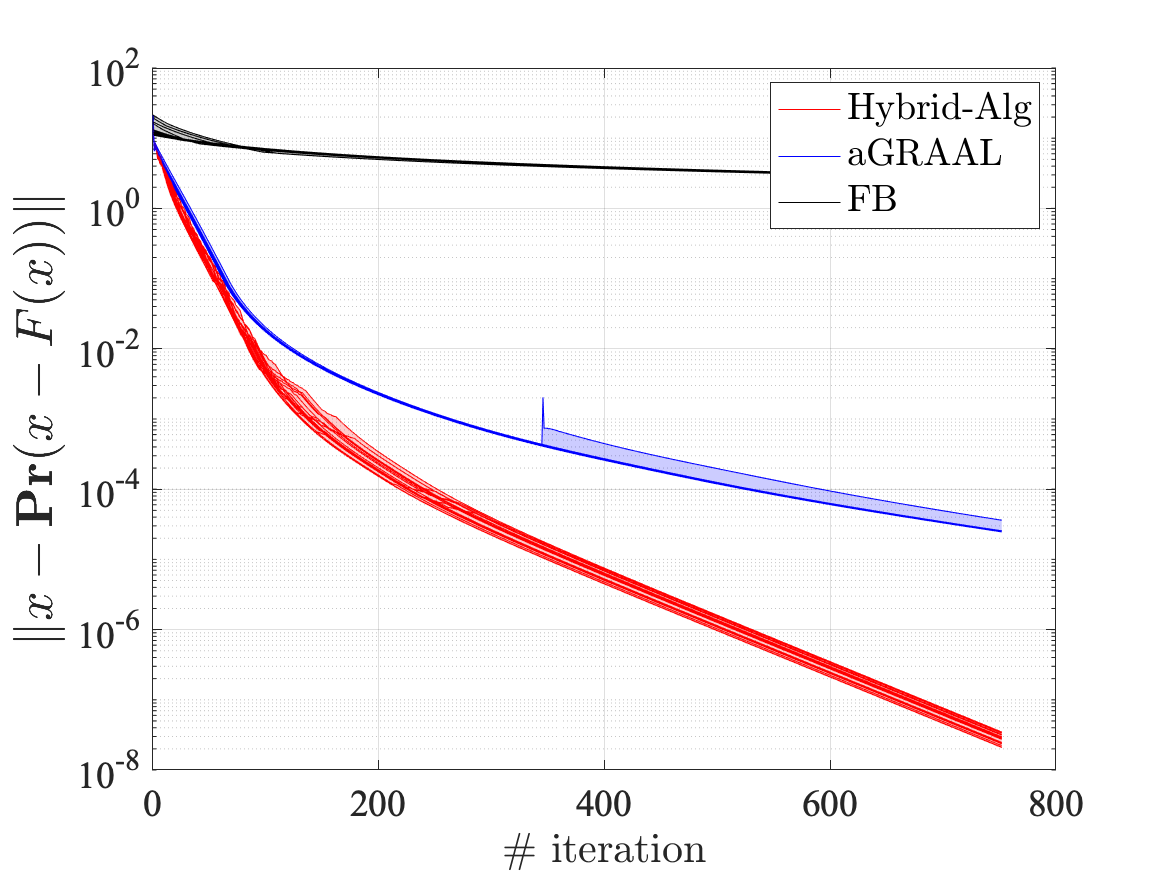}}
    \hspace{2mm}
    \subfloat[$K_i = \text{randi}\left((80, 120),1,n\right)$.]{\label{fig-conv-smp:3}\includegraphics[scale=0.35]{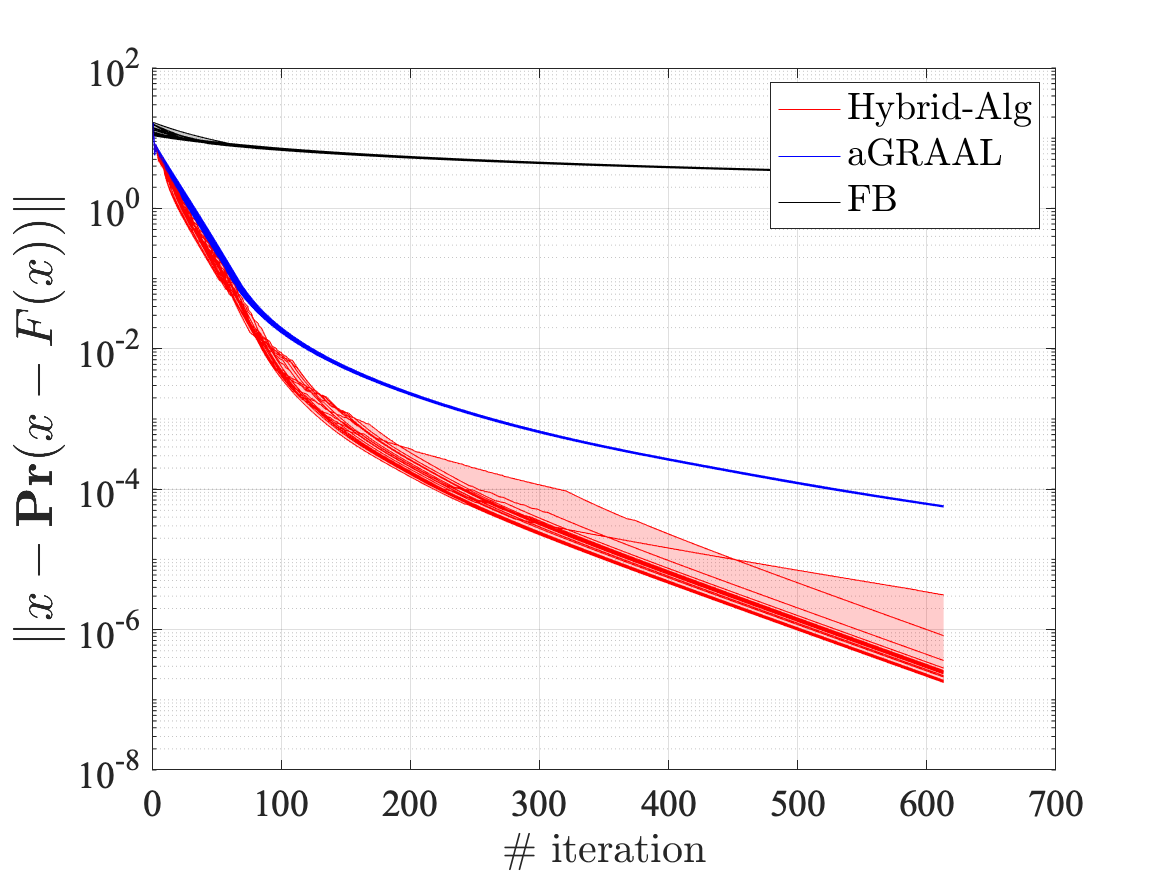}}
    \caption{Residual of variational inequality problem with different number of samples \SG{$K_i$ per agent} and radii $\epsilon_i$ chosen according to the distribution $\varepsilon \cdot U[1, 5]$ with fixed radius $\varepsilon = 0.01$.}
    \label{fig-conv-smp}
\end{figure}
\subsection{Numerical Example}

In this section, we reformulate a case study of the distributionally robust game in (\ref{DR_formulation}), under Assumption \ref{function_forms}, as a variational inequality problem and solve it using \RB{FB (with $\tau = 0.001$)}, aGRAAL, and Hybrid-Alg. The key difference between \RB{aGRAAL and Hybrid-Alg} algorithms is that, unlike aGRAAL, which uses a fixed momentum parameter, Hybrid-Alg employs a variable momentum parameter. We believe this is a testimony to the potential of switching the momentum parameter between a small (used in aGRAAL) and a large value, which has a significant impact on convergence speed. In particular, having larger, variable momentum parameter in Algorithm \ref{{Hybrid-Alg}2} makes $\bar{x}_k$ closer to the most recent iterate, $x_k$, rather than $\bar{x}_{k-1}$, which allows us to estimate the local Lipschitz constant of the corresponding VI mapping $F$ more precisely compared to aGRAAL.

\RB{Briefly speaking, the following two equations are evaluated in each iteration of Hybrid-Alg (Algorithm \ref{{Hybrid-Alg}2}) to assess sufficient decrease and to determine whether the large or the small momentum parameter should be used.}
\begin{align}\label{eq1:hyb-alg}
\frac{\theta_{k-1}}{2}\n{x^{k}-x^{k-1}}^2 &-\dfrac{\tau_k}{\tau_{k-1}}\phi_k \n{x^k-\x^k}^2 \nonumber + \bigl(\dfrac{\tau_k}{\tau_{k-1}}\phi_k - 1 - \frac{1}{\phi_{k+1}}\bigr) \n{x^{k+1}-\x^{k}}^2  \nonumber\\
&- \bigl(\dfrac{\tau_k}{\tau_{k-1}}\phi_k - \theta_k\bigr)\n{x^{k+1}-x^k}^2  - \frac{\theta_k}{2}\n{x^{k+1}-x^k}^2.
\end{align}
\begin{align}\label{eq2:hyb-alg}
\color{black} -\dfrac{\tau_k \phi_k}{\tau_{k-1}} \n{x^k-\x^k}^2 +\bigl(\dfrac{\tau_k\phi_k}{\tau_{k-1}} - 1 - \frac{1}{\phi_{k+1}}\bigr) \n{x^{k+1}-\x^{k}}^2 - \bigl(\dfrac{\tau_k\phi_k}{\tau_{k-1}} - \theta_k\bigr)\n{x^{k+1}-x^k}^2.
\end{align}
\par
\begin{figure}[t]
    \centering
    \subfloat[Different radii $\varepsilon$.]{\label{fig-val:1}\includegraphics[scale=0.35]{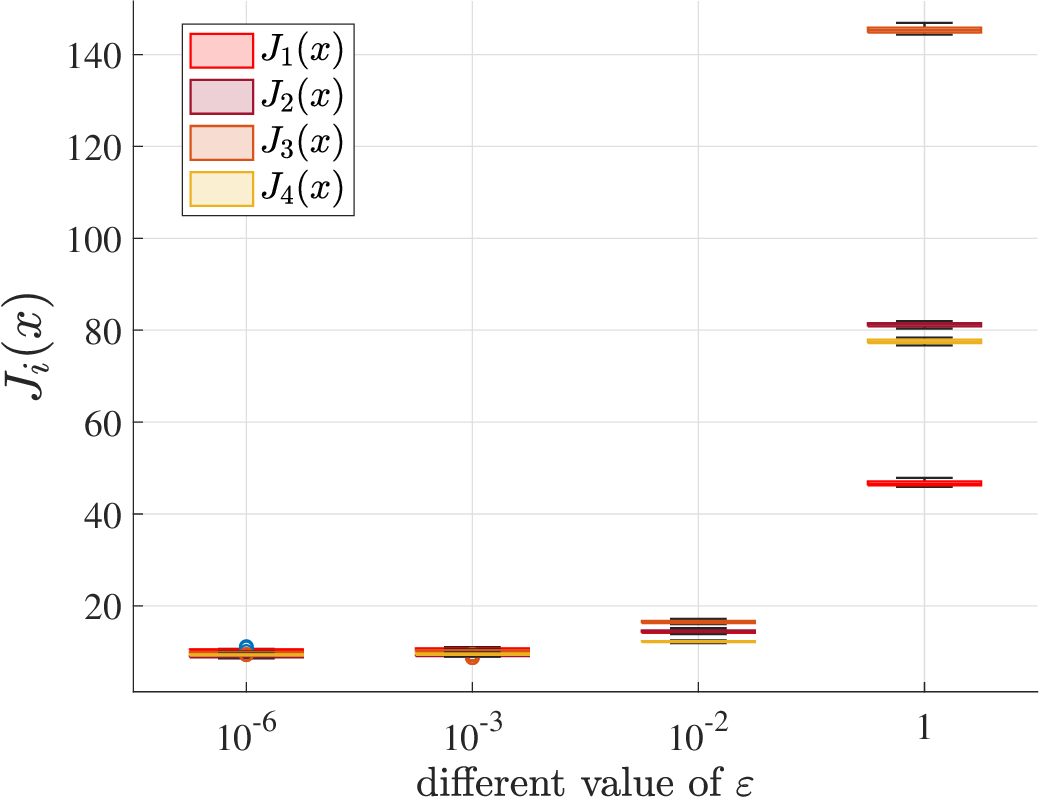}}
        \hspace{10mm}
    \subfloat[Different number of samples.]{\label{fig-val:2}\includegraphics[scale=0.35]{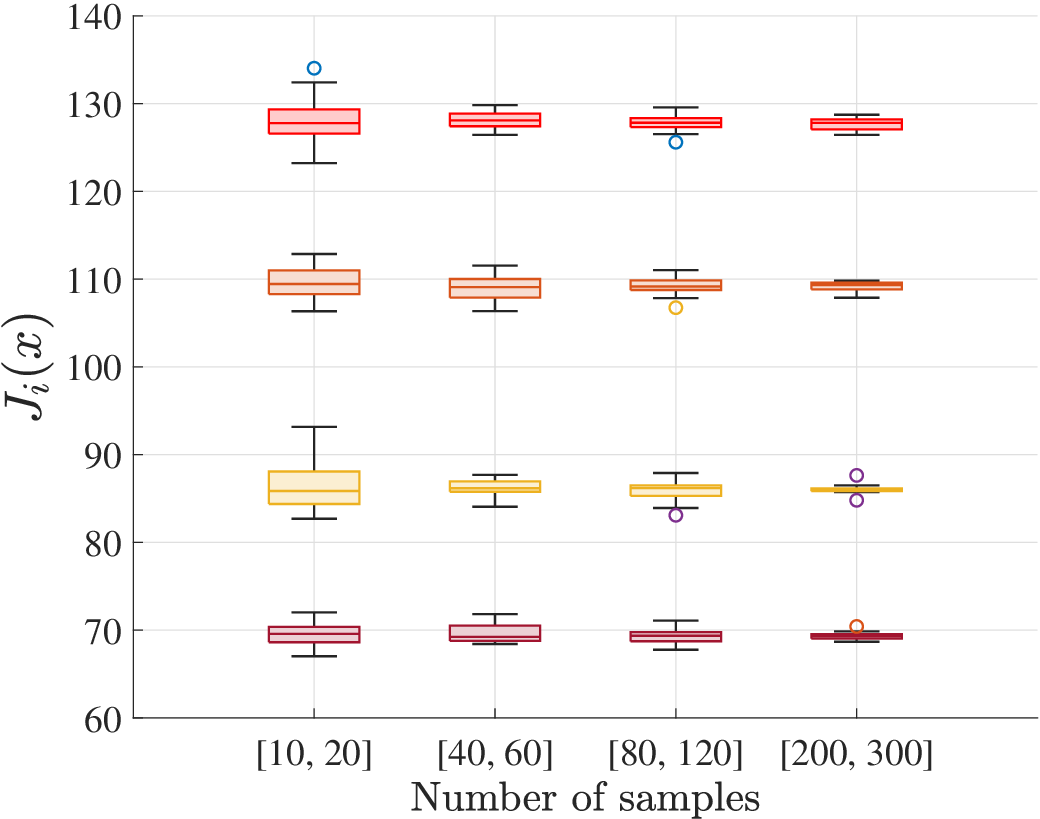}}
    \caption{Cost function values of four agents at the equilibrium for 10 different scenarios represented by box plots with (a) different radii 
    (b) Both different samples and different number of samples per agent.
    }
    \label{fig-conv-illus}
\end{figure}

 For the simulation, the parameters in the problem are generated as follows: 
Each drawn sample \( \xi_i^{(k_i)} \) is generated from the uniform distribution with support set [0,1], while  $P_i$ is given by $P_i(x)=\sum_{i \in \pazocal{N}}a_ix_i$.  The values  $a_i$, and the eigenvalues of $D_i$ in the reformulation (\ref{reformulation_VI}) are randomized. Each agent's Wasserstein radius $\varepsilon_i$ is chosen randomly according to the distribution $\varepsilon \cdot U[1,5]$, where $\varepsilon$ takes fixed values in $\{10^{-6}, 10^{-3}, 10^{-2}, 1\}$ and $U[1,5]$ is a uniform discrete distribution with support set \{1,2, \dots, 5\}. Figure \ref{fig-conv-radii} shows the residual of the corresponding mapping $F$ for the illustrative example for different Wasserstein radii and a fixed number of samples. We note that the convergence rate of both algorithms is almost linear, which illustrates that, even though the VI mapping can be nonmonotone, fast solutions can be obtained using both algorithms. Figure  \ref{fig-conv-smp} shows the residual for an increasing number of invividual data for each agent and individual radii per agent. The number of each agent's samples for each case study is drawn from a discrete integer distribution in $[10, 20]$, $[40, 60]$ and $[80, 120]$, respectively.  Note that even if we increase the number of samples, the convergence rate does not change, thus leading to results that scale well with the sample size. Finally, Figure \ref{fig-conv-illus} illustrates how the cost of each agent at the equilibrium is affected by the Wasserstein radii and the number of samples of each agent for 10 different problem instances represented by boxplots. In Figure \ref{fig-conv-illus}(a) we consider different radii $\varepsilon_i$ per agent obtained from the distribution $\varepsilon\cdot  U[1,5]$, where $\varepsilon\in \{10^{-6}, 10^{-3}, 10^{-2}, 1\}$ to investigate the effect of increasing Wasserstein radii on the cost of each agent; In  Figure \ref{fig-conv-illus}(b) the number of samples per agent follows a uniform distribution with support sets $\{[10,20], [40, 50], [80, 120], [200, 300]\}$ per case study to investigate the effect an increasing number of samples has on the cost of each agent. \par 
We observe that as we increase the value of the radii, the cost functions of each agent are higher representing a more conservative but robust behaviour against distrubutional shifts. Finally, for fixed radii, as the number of samples increases, the empirical variance of the costs decreases as well, as a result of a more accurate estimation of the probability distribution, used as \SG{the center of each ambiguity set}. 

\subsection{Risk-Aware Portfolio Allocation under Market Uncertainties and Behavioural Influences}
\begin{figure}[h]
    \centering
    \subfloat[$\varepsilon = 10^{-6}$.]{\label{fig-port-conv:1}\includegraphics[scale=0.35]{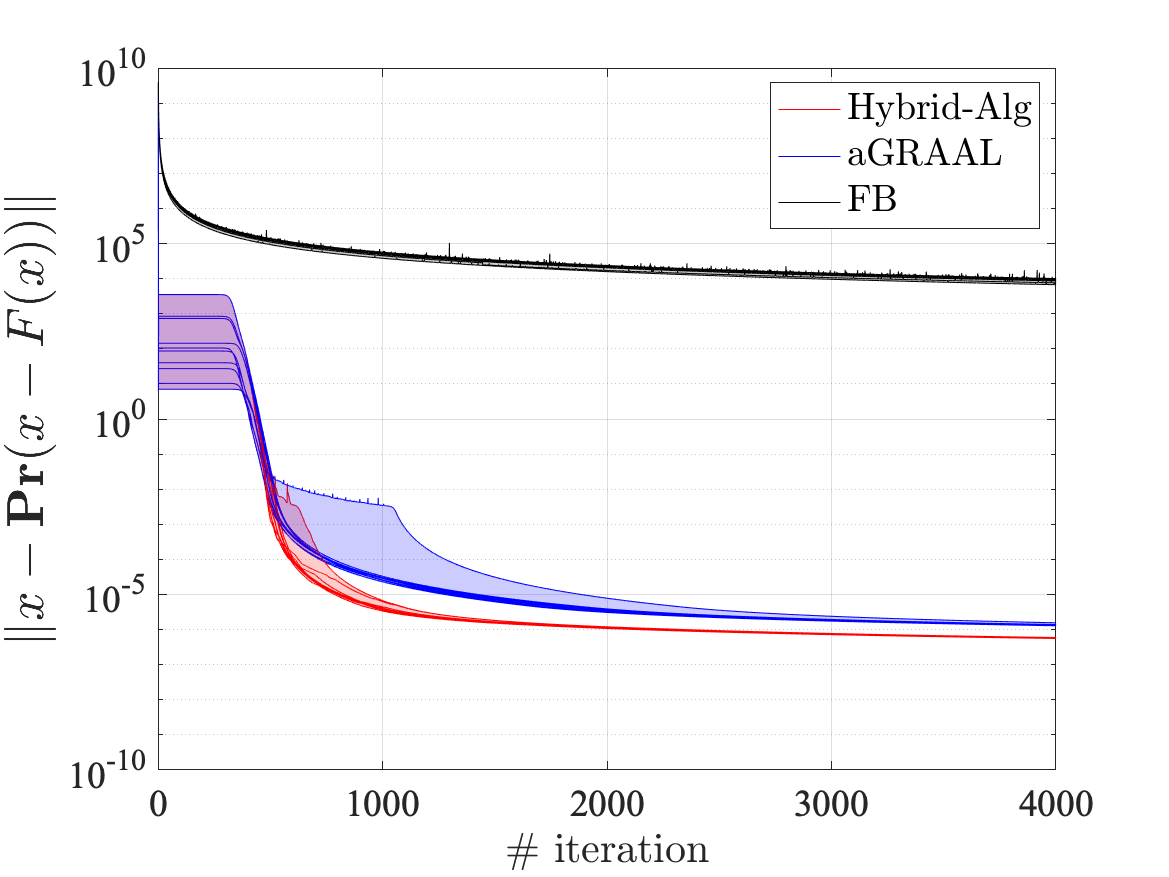}}
       \hspace{2mm}
    \subfloat[$\varepsilon = 10^{-2}$.]{\label{fig-port-conv:2}\includegraphics[scale=0.35]{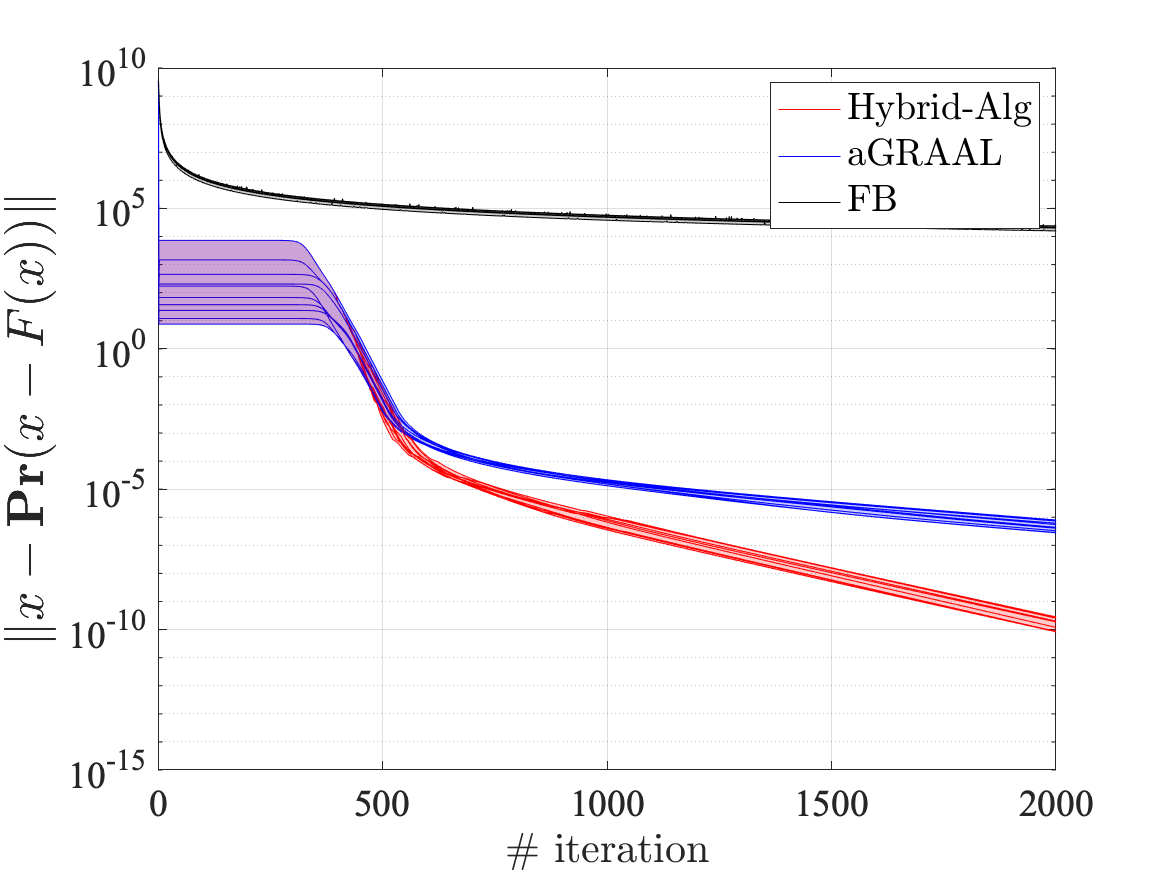}}
    \hspace{2mm}
    \subfloat[$\varepsilon = 1$.]{\label{fig-port-conv:3}\includegraphics[scale=0.35]{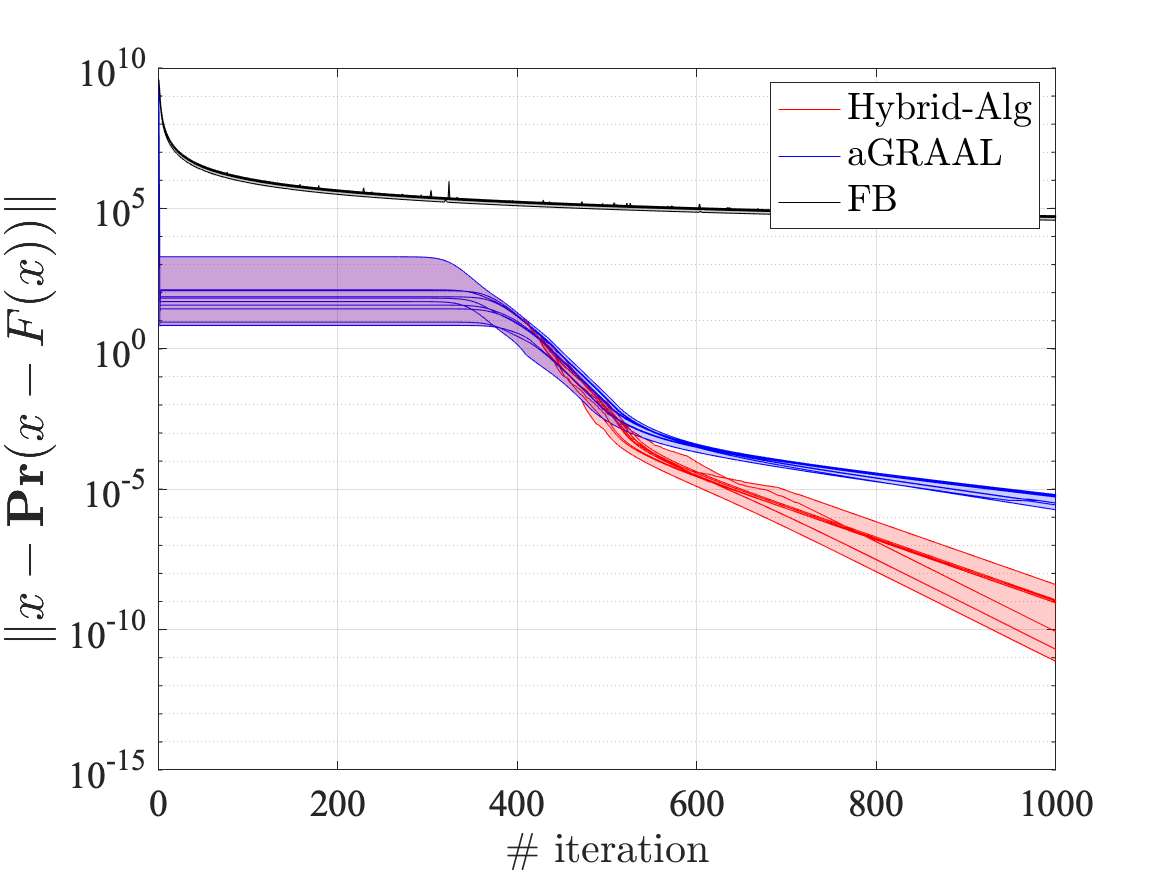}}
    \caption{Residual for different instances of the portofolio allocation problem with each agent's $\varepsilon_i$ chosen according to $\varepsilon \cdot U[1, 5]$ with $\varepsilon \in \{10^{-6}, 10^{-2}, 1\}$.}
    \label{fig-portfilio-conv}
\end{figure}
\begin{figure}[h]
    \centering
    \includegraphics[scale=0.35]{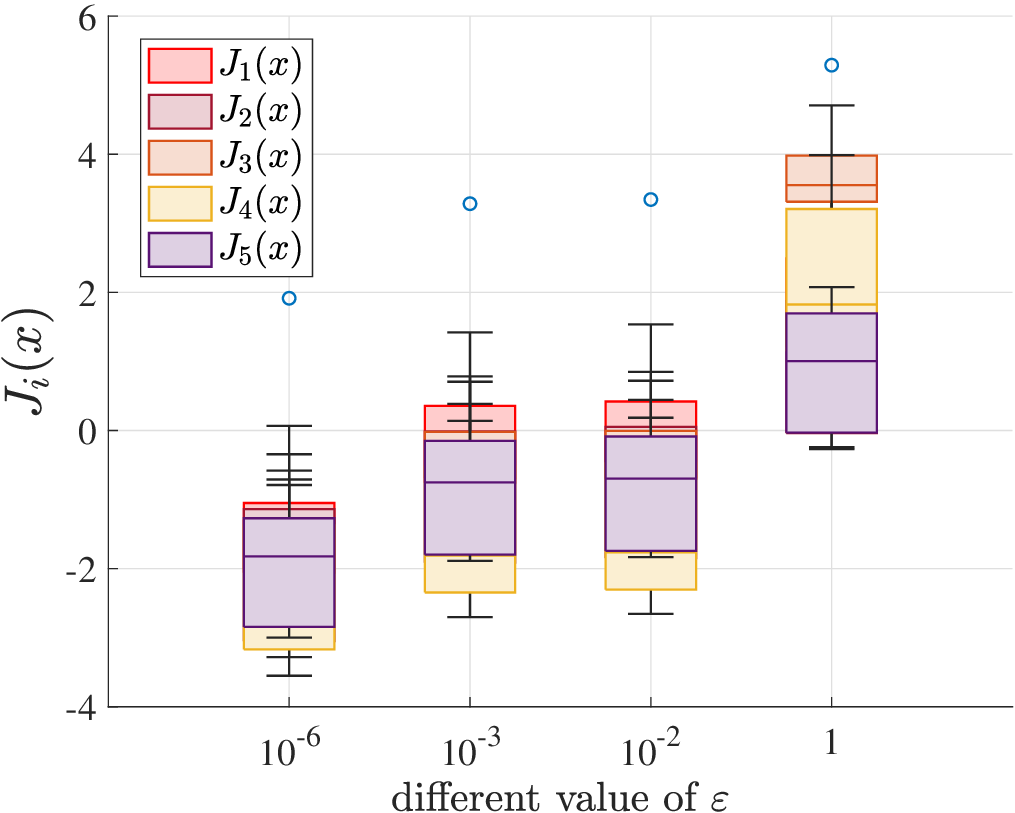}
    \caption{Cost function values of four agents at the equilibrium for 10 different scenarios represented by box plots for different radii $\varepsilon_i$ per agent obtained from the uniform distribution $\varepsilon \cdot U[1,5]$, where $\varepsilon\in \{10^{-6}, 10^{-3}, 10^{-2}, 1\}$ to investigate the effect of increasing Wasserstein radii on the cost of each agent.}
    \label{fig-conv}
\end{figure}
We consider a multi-investor robust portfolio allocation problem, where each investor \(i \in \mathcal{N}\) allocates capital seeking to maximize their profits or minimize their costs taking into account their exposure to market risks. The decision variable for each investor is their portfolio allocation \(x_i \in X_i\), where \(X_i\) represents the set of feasible portfolios for investor \(i\), \GP{normalized to a simplex representing the percentage of capital split among investments}. \GP{Furthermore, we wish to model behavioural impacts of other investors onto each individual investor}. \GP{Finally, we consider that agents are not only aware of the possible high variance of market uncertainties, but also aware that, when multiple investors accumulate to a single asset, this could lead to market bubbles which affects the returns from such investments.} Thus, each investor’s objective, given the other investors' strategies $x_{-i}$, is defined according to the following optimization problem:
\begin{align}
  \quad \min_{x_i \in X_i} \max_{\mathbb{Q}_i \in \mathcal{P}_i} \left\{ x_i^\top C_{ii}x_i + x_i^\top C_{ij}x_j - r^\top_ix_i  +\mathbb{E}_{\xi_i \sim \mathbb{Q}_i} \left[ \xi_i^\top Q_i \xi_i + P_i(x) \xi_i \right] \right\}.\nonumber
\end{align}
The term $r^\top_ix_i$ represents the deterministic part of the returns based on the allocation of capital to assets. The quadratic deterministic terms models (possible) behavioural coupling due to competition of the investors according to performance metrics often used to make such investments. 
The ambiguity set \( \mathcal{P}_i \) models investor \(i\)’s ambiguity in the distribution of uncertain market parameters affecting the returns. 
 The term \( \xi_i^\top Q_i \xi_i \) represents \(i\)’s aversion to volatility, indicating each agent's individual sensitivity to uncertain fluctuations. The term \( P_i(x) \xi_i \), where \( P_i(x) = \sum_{j \in \mathcal{N}} x_j \), models herding behavior, where multiple investors investing heavily in the same assets increase asset-specific risks. This crowding effect can drive prices up, raising the risk of market bubbles. In Figure \ref{fig-portfilio-conv}, we set the Wasserstein radii at $\varepsilon\in \{10^{-6}, 0.01, 1\}$ and consider 10 different instances of the problem with different values of matrices $Q_i, C_{ii}, C_{ij}, r_i$ and different multi-samples per agent $i \in \mathcal{N}$ obtained from different $t$-distributions. Note that even though the mapping is in general nonmonotone, most case studies lead to satisfactory (mostly linear) convergence results with both schemes, Hybrid-Alg (red lines) and aGRAAL (blue lines). In most of the case studies, the superiority of the hybrid algorithm is evident.
Figure \ref{fig-conv} shows the values of the cost functions of the agents at the equilibrium point for those 10 different instances represented by a box plot. Even though the problem is nonmonotone, increasing the Wasserstein radius of the agents leads in general to a larger value of the cost function. 
\section{Conclusion}

\GP{This work explores data-driven distributionally robust games using individual Wasserstein ambiguity sets and private data, thus allowing agents to develop their own personalized risk-averse decisions. We reformulate a seemingly-infinite dimensional game into a data-driven finite-dimensional variational inequality problem, which  evidently enjoys data-scalability properties. Future work will focus on introducing coupling constraints to our model. Extending on that we wish to investigate this problem under the presence of distributionally robust chance constraints coupling the agents decisions and in particular, whether certain assumptions such as linearity of the constraints can aid in obtaining a satisfactory reformulation or approximation of the original game. }

\begin{acknowledgements}
\GP{This research is partially supported  by the ERC under project COSMOS (802348).} 
\end{acknowledgements}
\bibliographystyle{plain}
\bibliography{biblio_coop.bib}
\end{document}